\documentclass[11pt]{article}

\usepackage[utf8]{inputenc}

\usepackage[T1]{fontenc} 
\usepackage{soulutf8} 
\sethlcolor{yellow} 

\usepackage[a4paper, margin=2.7cm]{geometry}
\usepackage[nottoc]{tocbibind}
\usepackage{mathtools, amsthm, amsfonts, amssymb} 
\usepackage{mathrsfs}           
\usepackage{enumitem} 
\usepackage{titlesec}
\usepackage[colorlinks=true,linkcolor=blue,citecolor=blue,urlcolor=blue,breaklinks]{hyperref}

\usepackage{graphicx} 
\usepackage{mathtools}
\usepackage{mathabx}
\usepackage{esint} 
\usepackage{bbm}
\usepackage{bigints}
\usepackage[dvipsnames]{xcolor}
\usepackage[toc,page]{appendix}
\usepackage{comment} 

\allowdisplaybreaks 

\usepackage{multirow}

\usepackage{tikz}
\usetikzlibrary{decorations.pathreplacing}

\numberwithin{equation}{section}

\DeclareMathOperator*{\fiint}{\ensuremath{\iint\text{\kern-1.36em{\raisebox{5.87pt}{\rotatebox{-93}{$\setminus$}}}}}}

\DeclareMathOperator*{\supess}{\ensuremath{ess{\text{ }}sup}}

\newtheorem{theorem}{Theorem}[section]
\newtheorem{lemma}[theorem]{Lemma}

\newtheorem{proposition}[theorem]{Proposition}
\newtheorem{assumption}[theorem]{Assumption}

\newtheorem{definition}[theorem]{Definition}

\newcommand{\Tr}{{\mathrm{Tr}}}
\newcommand{\D}{{\mathcal{D}}}

\newcommand{\rank}{\mathrm{Rank}}

\newtheorem{remark}[theorem]{Remark}

\providecommand\given{}
   \newcommand\SetSymbol[1][]{%
      \nonscript\:#1\vert
      \allowbreak
      \nonscript\:
      \mathopen{}}
   \DeclarePairedDelimiterX\Set[1]\{\}{%
      \renewcommand\given{\SetSymbol[\delimsize]}
      #1
}

\newcommand{\interior}[1]{%
  {\kern0pt#1}^{\mathrm{o}}%
}

\usepackage[style=numeric-comp,maxbibnames=9,natbib=true]{biblatex}
\AtEveryBibitem{%
  \clearfield{issn} 
  \clearfield{doi} 

  \ifentrytype{online}{}{
    \clearfield{url}
  }
}

\renewbibmacro*{name:andothers}{
  \ifboolexpr{
    test {\ifnumequal{\value{listcount}}{\value{liststop}}}
    and
    test \ifmorenames
  }
    {\ifnumgreater{\value{liststop}}{1}
       {\finalandcomma}
       {}%
     \andothersdelim\bibstring[\emph]{andothers}}
    {}}

\usepackage{hyperref}
\newlist{primenumerate}{enumerate}{1}
\setlist[primenumerate,1]{label={(H\arabic*$'$})}
\title{Properties of periodic Dirac--Fock functional and minimizers}
\author{
Isabelle Catto\footnote{ \textsc{Isabelle Catto, CEREMADE, Université Paris-Dauphine, Université PSL, CNRS, 75016 Paris, France}
  \textit{E-mail address}: \texttt{\href{mailto:catto@ceremade.dauphine.fr}{catto@ceremade.dauphine.fr}}}
\and Long Meng\footnote{ \textsc{Long Meng, CERMICS, \'Ecole des ponts ParisTech, 6 and 8 av. Pascal, 77455 Marne-la-Vall\'ee, France}
  \textit{E-mail address}:\texttt{\href{long.meng@enpc.fr}{long.meng@enpc.fr}}}
  }
\date{}

\begin{document}

\maketitle

\begin{abstract}
Existence of minimizers for the Dirac--Fock model in crystals  was recently proved by Paturel and S\'er\'e and the authors \cite{crystals} by a retraction technique due to S\'er\'e \cite{Ser09}. In this paper, inspired by Ghimenti and Lewin's result \cite{ghimenti2009properties} for the periodic Hartree--Fock model, we prove that the Fermi level of any periodic Dirac--Fock minimizer is either empty or totally filled when $\frac{\alpha}{c}\leq C_{\rm cri}$ and $\alpha>0$. Here $c$ is the speed of light, $\alpha$ is the fine structure constant, and $C_{\rm cri}$ is a constant only depending on the number of electrons and on the charge of nuclei per cell. More importantly, we provide an explicit upper bound for $C_{\rm cri}$.

Our result implies that any minimizer of the periodic Dirac--Fock model is a projector when $\frac{\alpha}{c}\leq C_{\rm cri}$ and $\alpha>0$. In particular, the non-relativistic regime (i.e., $c\gg1$) and the weak coupling regime (i.e., $0<\alpha\ll1$) are covered. 
   
The proof is based on a delicate study of a second-order expansion of the periodic Dirac--Fock functional composed with the retraction used in \cite{crystals}.
\end{abstract}

\section{Introduction}
This paper is devoted to the study of the properties of the Dirac--Fock ground-state in  perfect crystals; that is,  in the periodic case.

The Hartree--Fock (HF) model is commonly used in non-relativistic chemistry and quantum physics to calculate ground- or bound-state energies of atoms and molecules. In this model, the state of the electrons is represented by a so-called density matrix  $\gamma$ which is a self-adjoint trace-class operator $0\leq \gamma\leq 1$ acting on the space $L^2(\mathbb{R}^3;\mathbb{C})$. Its finite trace represents the number  $N $ of electrons ($N\in \mathbb{N}^+$). When the nuclear charge $Z>N-1$, existence of ground-states for this model expressed in terms of $N$-particle wave-functions goes back to Lieb and Simon \cite{LieSim-77}, and has been extended later to excited states by Lions~\cite{Lions-87}, (see also the recent review paper by Bach \cite{Bach2022} and the references therein). Existence of minimizers for the HF functional involving one-particle density matrix is due to Lieb~\cite{lieb1981variational} (see also Bach and al.~\cite{bach1992error}). Additionally, it is shown in \cite{lieb1981variational}  that any HF minimizer $\gamma$ is automatically a projector of the form $\gamma=\sum_{n=1}^N\left|\psi_n\right>\left<\psi_n\right|$ with the $\psi_i$'s being the  eigenfunctions  corresponding to the smallest eigenvalues (counted with multiplicity) of the mean-field self-adjoint HF operator $H_{\gamma}$,
\[
H_{\gamma}\psi_n=\epsilon_n\psi_n,
\]
with spectrum $\sigma(H_\gamma)=\{\epsilon_1\leq \epsilon_2\leq\cdots\leq \epsilon_N\leq \cdots\}\bigcup [0,+\infty)$ where $\epsilon_j<0$ for any $j\in \mathbb{N}^+$. Furthermore, Bach, Lieb, Loss and Solovej proved that shells are always completely filled in the HF model~\cite{bach1994unfill}. Mathematically, this property writes  $\epsilon_N<\epsilon_{N+1}$. In particular, any minimizer of the HF functional solves the following self-consistent equation
\[
\gamma=\mathbbm{1}_{(-\infty,\epsilon_{N}]}(H_{\gamma}),\]
where $\mathbbm{1}_I(H)$ denotes the spectral projection on the set $I\subset \mathbb{R}$ of the self-adjoint operator $H$. All these properties are important for the efficient numerical methods in the HF theory (see e.g., \cite{cances2000convergence}).

This result has been  later extended by Ghimenti and Lewin in \cite{ghimenti2009properties} to the periodic HF model for neutral crystals introduced and studied  by Catto, Le Bris and Lions in \cite{catto2001thermodynamic}. They proved that any minimizer $\gamma$ of the periodic HF energy is indeed always a projector (of infinite rank), that solves the self-consistent equation 
\begin{align}\label{eq:gamma-HF}
    \gamma=\mathbbm{1}_{(-\infty,\nu)(H_{\gamma})}+\epsilon\mathbbm{1}_{\{\nu\}}(H_{\gamma})
\end{align}
with $\epsilon\in \{0,1\}$, and $\nu$ may be an eigenvalue of the periodic HF operator $H_{\gamma}$ (with infinite multiplicity  due to the invariance by translations of the lattice.) In \cite{cances2008new}, a similar result was proved {by Canc\`es, Deleurence and Lewin} for the reduced HF model for crystals where the exchange term is neglected. Based on \cite{thomas1973time}, they show that the spectrum of the corresponding self-adjoint  operator is purely absolutely continuous. Hence $\nu$ cannot be an eigenvalue and one can take $\epsilon=1$ in \eqref{eq:gamma-HF}. In particular, there are no unfilled shells in the reduced HF theory. Unfortunately,  we do not know whether the spectrum of the periodic HF operator is also purely absolutely continuous because of the non-local feature of the exchange term.

 When heavy nuclei are involved (that is, $Z$ large), it is expected that the electrons closest to the nucleus move at very high velocities, thus requiring a relativistic treatment. It is widely believed that Quantum Electrodynamics (QED) is an adequate framework to deal with relativistic effects. It is shown that the shells are always completely filled for the Bogoliubov--Dirac--Fock model in QED \cite{hainzl2009existence}. The proof relies on the fact that the corresponding functional is bounded from below and on the positive definiteness of the nonlinear term. However, this theory leads to divergence problems: It is not easy to give meaning to different physical quantities appearing in QED, such as the energy and the  charge density of the vacuum. 

Alternatively, the Dirac--Fock model  (DF) for atoms and molecules is one of the most attractive models in relativistic computational chemistry. It is a variant of the HF model in which the Laplace operator $-\frac{1}{2}\Delta$ entering the kinetic energy term is replaced by the free Dirac operator $\D$. Unlike the QED models, the DF functional is not bounded from below. It is therefore difficult to give a rigorous definition of the ground-state energy. However, existence results of critical points -- that is,  solutions to the DF equations -- can be found in \cite{esteban1999solutions,paturel2000solutions}. These solutions provide an infinite number of finite rank projectors as critical points of the DF functional. It is also proved in \cite{esteban2001nonrelativistic} that, up to subsequences, the projector with the smallest energy among these critical points converges to a minimizer of the Hartree-Fock energy in the non-relativistic limit; that is when the speed of light goes to infinity.

Recently, in the spirit of Lieb's variational principle (see, e.g., \cite{lieb1981variational,bach1992error}), S\'er\'e redefined the DF ground-state energy for atoms and molecules by using the density matrix formalism~\cite{Ser09}. Using a retraction technique, he proved that the DF ground-state energy admits a minimizer $\gamma$ on a suitable subset of density matrix, and that $\gamma$ satisfies the self-consistent equation
\begin{equation*}
    \gamma=\mathbbm{1}_{(0,\nu)}(\D_\gamma)+\delta,\quad \textrm{ with } 0\leq\delta\leq \mathbbm{1}_{\{\nu\}}(\D_{\gamma}),
\end{equation*}
for some Lagrange multiplier $\nu$ associated with the charge constraint. Later on, by using Séré's retraction technique, Meng~\cite{meng} justifies mathematically Mittleman's approach to the DF model: The DF model is an approximation of a max-min problem coming from the electron-positron field (see, e.g.,  \cite{barbaroux2005hartree,barbaroux2005some,huber2007solutions}). As a byproduct, he shows that the shells in the DF theory of  atoms and molecules are completely filled when the fine structure constant $\alpha$ is small enough or the speed of light $c$ is large enough under some conditions on $Z$ and $N$. This is an immediate consequence of a second-order expansion of a new DF functional: a composition of the DF functional with the retraction.

Finally,  one can construct the periodic DF model for crystals by replacing the Schr\"ondinger operator in the periodic HF model by the Dirac operator. Together with Paturel and S\'e\'e, we recently studied this new model in~\cite{crystals}. We show that the energy of the periodic DF model admits a minimizer that solves the self-consistent equation
\begin{equation*}
    \gamma=\mathbbm{1}_{(0,\nu)}(\D_\gamma)+\delta,\quad \textrm{ with } 0\leq\delta\leq \mathbbm{1}_{\{\nu\}}(\D_{\gamma}).
\end{equation*}
Here $\gamma$ and $\D_{\gamma}$ are the periodic density matrix and the periodic DF operator respectively. In the present paper, inspired by the results of one of us~\cite{meng}, we investigate the properties of the DF minimizers in crystals. We  mimic the proof of Ghimenti and Lewin for the periodic HF model~\cite{ghimenti2009properties} and obtain a  similar result. More precisely, we show that when $\alpha>0$ and $\frac{\alpha}{c}\leq C_{\rm cri}$ (with $C_{\rm cri}$ given in Lemma \ref{lem:contra2}), any minimizer $\gamma$ of the periodic DF ground-state energy is always a projector, and that it solves a self-consistent  equation of the form
\begin{align*}
    \gamma=\mathbbm{1}_{(0,\nu)}(\D_{\gamma})+\epsilon\mathbbm{1}_{\{\nu\}}(\D_{\gamma})
\end{align*}
with $\epsilon\in \{0,1\}$. The proof of Ghimenti and Lewin  is  based on the  local convexity of the periodic HF  functional. Since this property does not hold any longer here, we rather rely on a careful study of the second-order expansion of the periodic DF functional due to one of us~ \cite{meng}. 
\section{Description of the periodic DF model and main results}
The paragraph below is copied from \cite{crystals} for the reader's convenience. For the sake of simplicity, we only  consider the case of a cubic crystal with a single point-like nucleus per unit cell, that is located at the center of the cell. The reader should however keep in mind that the general case could be handled as well. Let $\ell>0$ denote the length of the elementary cell $Q_\ell=(-\frac{\ell}{2},\frac{\ell}{2}]^3$. The nuclei with positive charge $z$  are treated as classical particles with infinite mass that are located at each point of the lattice $\ell\,\mathbb{Z}^3$. The electrons are treated quantum mechanically through a periodic density matrix. The electronic density is modeled by a $Q_\ell$-periodic function whose $L^1$-norm over the elementary cell equals the ``number of electrons'' $q\in \mathbb{N}^+$ per cell -- the electronic charge per cell being equal to $-q$. In this paper, we set $q\in \mathbb{N}^+$ and $z\in \mathbb{R}^+$. 

In this periodic setting, the $Q_\ell$-periodic Coulomb  potential $G_\ell$ resulting from a distribution of point particles of charge $1$ that are periodically located at the centers of the cubic cells of the lattice is defined, up to a constant, by
\begin{equation}\label{eq:def-G}
    -\Delta G_\ell=4\pi\left[-\frac{1}{\ell^3}+\sum_{k\in\mathbb{Z}^3}\delta_{\ell k}\right].
\end{equation}
By convention, we choose $G_\ell$ such that
\begin{equation}\label{eq:constante-G}
    \int_{Q_\ell}G_\ell\,dx=0.
\end{equation}
The Fourier series of $G_\ell$ writes
\begin{equation}\label{3.eq:2.6}
    G_\ell(x)=\frac{1}{\pi \ell}\sum_{p\in\mathbb{Z}^3\setminus\{0\}}\frac{e^{\frac{2i\pi}{\ell}p\cdot x}}{|p|^2}, \quad\text {for  every } x\in \mathbb{R}^3.
\end{equation}

The free Dirac operator is defined by $
\D=-ic\sum_{k=1}^3\alpha_k\partial_k+c^2\beta$, 
with $4\times4$ complex matrix $\alpha_1,\alpha_2,\alpha_3$ and $\beta$, whose standard forms are $
\beta=\begin{pmatrix} 
\mathbbm{1}_2 & 0 \\
0 & -\mathbbm{1}_2 
\end{pmatrix}$, 
$\alpha_k=\begin{pmatrix}
0&\sigma_k\\
\sigma_k&0
\end{pmatrix}$ 
where $\mathbbm{1}_2$ is the $2\times 2$ identity matrix and the $\sigma_k$'s, for $k\in\{1,2,3\}$, are the well-known $2\times2$ Pauli matrix $
\sigma_1=\begin{pmatrix}
0&1\\
1&0
\end{pmatrix},\,
\sigma_2=\begin{pmatrix}
0&-i\\
i&0
\end{pmatrix}$, 
$\sigma_3=\begin{pmatrix}
1&0\\
0&-1
\end{pmatrix}.$ Here $c >0$ denotes the speed of light.

The operator $\D$ acts on $4-$spinors; that is, on functions from $\mathbb{R}^3$ to $\mathbb{C}^4$. It is self-adjoint on $L^2(\mathbb{R}^3;\mathbb{C}^4)$, with domain $H^1(\mathbb{R}^3;
\mathbb{C}^4)$ and form domain $H^{1/2}(\mathbb{R}^3;\mathbb{C}^4)$ (denoted by $L^2$, $H^1$ and $H^{1/2}$ in the following, when there is no ambiguity). Its spectrum is $\sigma(\D)=(-\infty,-c^2]\bigcup[+c^2,+\infty)$. Following the notation in \cite{esteban1999solutions,paturel2000solutions}, we denote by $\Lambda^+$ and $\Lambda^-=\mathbbm{1}_{L^2}-\Lambda^+$ respectively the two orthogonal projectors on $L^2(\mathbb{R}^3;\mathbb{C}^4)$ corresponding to the positive and negative eigenspaces of $\D$; that is
\[
\begin{cases}
\D\Lambda^+=\Lambda^+\D=\Lambda^+\sqrt{c^4-c^2\Delta}=\sqrt{c^4-c^2\Delta}\,\Lambda^+,\\
\D\Lambda^-=\Lambda^-\D=-\Lambda^-\sqrt{c^4-c^2\Delta}=-\sqrt{c^4-c^2\Delta}\,\Lambda^-.
\end{cases}
\]
According to the Floquet theory \cite{reed1978b}, the underlying Hilbert space $L^2(\mathbb{R}^3;\mathbb{C}^4)$ is unitarily equivalent to $L^2(Q_\ell^*)\bigotimes L^2(Q_\ell;\mathbb{C}^4)$, where $Q_\ell^*=[-\frac{\pi}{\ell},\frac{\pi}{\ell})^3$ is the reciprocal cell of the lattice, with volume $|Q_\ell^*|=(2\pi)^3/\ell^3$. (In the Physics literature $Q_\ell^*$ is known as the first Brillouin zone.) The Floquet unitary transform $\mathfrak{U}:\;L^2(\mathbb{R}^3;\mathbb{C}^4)\to L^2(Q_\ell^*)\bigotimes L^2(Q_\ell;\mathbb{C}^4) $ is given by
\begin{align}
    \mathfrak{U}:\;\phi \mapsto \fint_{Q_{\ell}^*} (\mathfrak{U}\phi)_\xi d\xi
\end{align}
with the shorthand $\fint_{\Omega}$ for $\frac{1}{|\Omega|}\int_\Omega$ and
\begin{equation}\label{eq:def_Bloch_transform}
(\mathfrak{U}\phi)_\xi=\sum_{k\in\mathbb{Z}^3}e^{-i\ell k\cdot \xi}\phi(\cdot+\ell\,k)
\end{equation}
for every $\xi\in Q^*_\ell$ and $\phi$ in $L^2(\mathbb{R}^3;\mathbb{C}^4)$. For every $\xi\in Q^*_\ell$, the function $(\mathfrak{U} \phi)_\xi$ belongs to the space
\[
L^2_\xi(Q_\ell;\mathbb{C}^4)=\Set*{\psi\in L^2_{\text{loc}}(\mathbb{R}^3 ; \mathbb{C}^4)\given e^{-i\xi\cdot x}\psi\,\,\text{is }Q_\ell\text{-periodic}},
\]
which will also be denoted by $L^2_\xi$ in the sequel. We write $L^2(\mathbb{R}^3;\mathbb{C}^4)=\fint_{Q_\ell^*}^\oplus L^2_\xi\,d\xi \cong L^2(Q_\ell^*)\otimes L^2(Q_\ell;\mathbb{C}^4)$ to refer to this direct integral decomposition of $L^2$ w.r.t. the Floquet transform $\mathfrak{U}$. Functions $\psi$ of this form are called Bloch waves or $Q_\ell$-quasi-periodic functions with quasi-momentum $\xi\in Q^*_\ell$. They satisfy 
\[
\psi(\cdot+\ell\,k)=e^{i\ell\,k\cdot\xi}\psi(\cdot), \text{ for every }k\in \mathbb{Z}^3
.\]
 In particular, when $\psi$ is $Q_\ell$-periodic (i.e., $\xi=0$), we denote
\[
L^2_{\rm per}(Q_\ell):=L^2_0(Q_\ell).
\]
The free Dirac operator can be rewritten accordingly as 
\begin{equation}\label{eq:somme_directe_Dirac}
    \D=\fint_{Q_\ell^*}^\oplus \D_\xi \,d\xi.
\end{equation}
where the $\D_\xi$'s are self-adjoint operators on $L^2_\xi(Q_\ell;\mathbb{C}^4)$ with respective domains $H_\xi^1(Q_\ell;\mathbb{C}^4)$ and respective form-domains $H_\xi^{1/2}(Q_\ell;\mathbb{C}^4)$.  Note that 
\begin{align*}
    \D_\xi^{\,2}=c^4-c^2\Delta_\xi,
\end{align*}
where $-\Delta=\fint_{Q_\ell^*}^\oplus-\Delta_\xi d\xi$. 

For almost every $\xi\in Q_\ell^*$,  the spectrum $\sigma(\D_\xi)$ of $\D_\xi$ is composed of two sequences of real eigenvalues $(d^-_n(\xi))_{n\geq 1}$ and $(d^+_n(\xi))_{n\geq 1}$ such that, for every $n\geq 1$,  
\begin{align*}
d^-_n(\xi)\leq -c^2,\quad \lim_{n\to\infty }d^-_n(\xi)=-\infty,\qquad d^+_n(\xi)\geq + c^2
,\quad \lim_{n\to\infty }d^+_n(\xi)=+\infty,
\end{align*}
and
\[
\sigma(\D)=\bigcup_{\xi\in Q_\ell^*}\sigma(\D_\xi)=\adjustlimits\bigcup_{\xi\in Q_\ell^*}\bigcup_{n\geq 1}\left\{d^-_n(\xi),d^+_n(\xi)\right\}=(-\infty,-c^2]\bigcup[+c^2,+\infty).
\]
\subsection{Functional framework}
As in \cite{crystals}, we now introduce various functional spaces for linear operators onto $L^2(Q_\ell;\mathbb{C}^4)$ and for operators onto $L^2(\mathbb{R}^3;\mathbb{C}^4)$ that commute with translations. 
Let $\mathcal{B}\left(E \right)$ 
be the set of bounded operators on a Banach space $E$ to itself. We use the shorthand $\mathcal{B}(L^2_\xi)$ for $\mathcal{B}(L^2_\xi(Q_\ell;\mathbb{C}^4))$. The space of bounded operators on $\fint_{Q_\ell^*}^\oplus L^2_\xi\,d\xi$ which commute with the translations of $\ell \mathbb{Z}^3$ is denoted by $Y$. It is isomorphic to  $L^\infty(Q_\ell^*;\mathcal{B}(L^2_\xi))$, and, for every  $h=\fint_{Q_\ell^*} h_\xi\,d\xi \in Y$, 
\[
\Vert h\Vert_{Y}=\supess_{\xi\in Q_\ell^*}\Vert h_\xi\Vert_{\mathcal{B}(L^2_\xi)}=\|h\|_{\mathcal{B}(L^2(\mathbb{R}^3;\mathbb{C}^4))}
\]
(see \cite[Theorem XIII.83]{reed1978b}). In this paper, we also use another norm on $Y$ which is defined by
\[
\|h\|_\mathcal{Y}=\sup_{\xi\in Q_\ell^*}\|h\|_{\mathcal{Y}(\xi)},
\]
with
\begin{align*}
 \|h\|_{\mathcal{Y}(\xi)}=\fint_{Q_{\ell}^*}\frac{\|h_{\xi'}\|_{\mathcal{B}(L^2_{\xi'})}}{|\xi-\xi'|^2}\,d\xi'.
\end{align*}
This convolution-type norm plays a critical role in this paper (see Lemma \ref{lem:contra2}). It is easy to see that, there exists a positive constant $C_\mathcal{Y}$ such that for every $h\in Y$,
\begin{align}\label{eq:T-Y}
    \|h\|_\mathcal{Y}\leq C_\mathcal{Y}\,\|h\|_Y.
\end{align}
In addition, we also use the rescaled $c$-dependent norms $\|h\|_{Y_c}:=c\|h\|_Y$ and $\|h\|_{\mathcal{Y}_c}:=c\|h\|_\mathcal{Y}$.

For $s\in [1,\+\infty)$ and $\xi\in Q_\ell^*$, we now define
 \[
\mathfrak{S}_{s}(\xi):=\Set*{h_\xi\in\mathcal{B}(L^2_\xi)\given \Tr_{L^2_\xi}(|h_\xi|^s)<\infty}
\]
endowed with the norm
\[
\|h_\xi\|_{\mathfrak{S}_s(\xi)}=\left(\Tr_{L^2_\xi}(|h_\xi|^s)\right)^{1/s}.
\]
We denote by $\mathfrak{S}_\infty(\xi)$ the subspace of compact operators in $\mathcal{B}(L^2_\xi)$, endowed with the operator norm  $\Vert\cdot\Vert_{\mathcal{B}(L^2_\xi)}$. Analogously, for $t\in [1,+\infty]$, we define
\[
\mathfrak{S}_{s,t}:=\Set*{ h=\fint^\oplus_{ Q^*_\ell} h_\xi \,d\xi \given h_\xi\in \mathfrak{S}_s(\xi)\text{ a.e., }\xi\in Q_\ell^*, \int_{Q_\ell^*}\|h_\xi\|_{\mathfrak{S}_s(\xi)}^t\,d\xi<\infty}
\]
endowed with the usual norm of $L^t(Q_\ell^*;\mathfrak{S}(\xi)^s)$:
\[
\|h\|_{\mathfrak{S}_{s,t}}=\left(\fint_{Q_\ell^*}\|h_\xi\|_{\mathfrak{S}_s(\xi)}^t d\xi\right)^{1/t}.
\]
We also define
\[
X_c^\alpha(\xi)=\Set*{h_\xi\in\mathcal{B}(L^2_\xi)\given |\D_\xi|^{\alpha/2}h_\xi|\D_\xi|^{\alpha/2}\in \mathfrak{S}_1(\xi)}
\]
endowed with the $c$-dependent norm
\[
\|h_\xi\|_{X_c^\alpha(\xi)}=\left\Vert|\D_\xi|^{\alpha/2}h_\xi|\D_\xi|^{\alpha/2}\right\Vert_{\mathfrak{S}_1(\xi)},
\]
and
\[
X^\alpha_{c,s}:=\Set*{h=\fint^\oplus_{ Q^*_\ell}  h_\xi \,d\xi \given h_\xi\in \mathfrak{S}_1(\xi) \text{ a.e. }\xi\in Q_\ell^*,  \int_{Q^*_\ell}\||\D_\xi|^{\alpha/2}h_\xi|\D_\xi|^{\alpha/2}\|_{\mathfrak{S}_1(\xi)}^s d\xi<\infty}
\]
endowed with the norm
\begin{equation*}
\|h\|_{X_{c,s}^\alpha}=\left(\fint_{Q_\ell^*}\||\D_\xi|^{\alpha/2}h_\xi|\D_\xi|^{\alpha/2}\|^s_{\mathfrak{S}_1(\xi)} d\xi\right)^{1/s}.
\end{equation*}
The space $X_c := X_{c,1}^1$ plays an important role in the definition of periodic density matrix, see Definition~\ref{def:one-part-density} below. For convenience, we use the notation $X_c(\xi)$ for $X^1_c(\xi)$. On $X_c$, we  will also use the norm (independent of $c$)
\[
\|\gamma\|_X:=\|(1-\Delta)^{1/4}\gamma(1-\Delta)^{1/4}\|_{\mathfrak{S}_{1,1}}.
\]
The norm on the intersection of any two functional spaces $A$ and $B$ will be defined by 
\[
\|\gamma\|_{A\bigcap B}=\max\{\|\gamma\|_{A},\|\gamma\|_{B}\}.
\]
We are now in the position to set the DF model for crystals. 
\subsection{The periodic  DF model}
We start with the following.
\begin{definition}[Periodic one-particle density matrix, \cite{crystals}]\label{def:one-part-density} We denote by $\mathcal{T}$ the set of $Q_\ell$-periodic one-particle density matrix
\[
\mathcal{T}:=\Set*{\gamma\in X_c\bigcap Y\given \gamma^*=\gamma,0\leq \gamma\leq\mathbbm{1}_{L^2(\mathbb{R}^3)}}.
\]
\end{definition}
We are particularly interested in the following subsets:
\[
\Gamma_{q}:=\Set*{\gamma\in \mathcal{T} \given \fint_{Q_\ell^*} \Tr_{L^2_\xi}(\gamma_{{*},\xi})\,d\xi= q}
\]
and 
\[
\Gamma_{\leq q}:=\Set*{\gamma\in \mathcal{T} \given  0\leq \fint_{Q_\ell^*} \Tr_{L^2_\xi}(\gamma_{{*},\xi})\,d\xi\leq q}.
\]
When $q$ is an integer, $ \Gamma_{q}$ (resp. $\Gamma_{\leq q}$) is the set of all DF states of a system of exactly $q$ (resp. at most $q$) electrons per unit cell. 

For $\gamma\in \Gamma_{\leq q}$, the periodic DF energy is defined by
\begin{align}
\mathcal{E}(\gamma)
 &=\fint\limits_{Q_\ell^*}\Tr_{L^2_\xi}[\D_\xi\gamma_\xi]\,d\xi- z\int\limits_{Q_\ell}G_\ell(x)\rho_\gamma(x)\,dx+\frac{\alpha}{2}\iint\limits_{Q_\ell\times Q_\ell}\rho_\gamma(x)G_\ell(x-y)\rho_\gamma(y)\,dxdy \notag\\
& \qquad \quad\qquad {} -\frac{\alpha}{2}\fiint_{Q_\ell^*\times Q_\ell^*}\,d\xi d\xi'\iint\limits_{Q_\ell\times Q_\ell}\Tr_{4}{[\gamma_\xi(x,y)\gamma_{\xi'}(y,x)]}W_\ell^\infty(\xi-\xi',x-y)\,dxdy.\label{3.energy}
\end{align}
In the  above definition of the energy, the so-called fine structure constant $\alpha$ is a dimensionless positive constant.

The potential $W_\ell^\infty$ that enters the definition of the last term, called the ``exchange term'', is defined by
\begin{equation}\label{3.eq:2.8}
W_\ell^\infty(\eta,x)=\sum_{k\in\mathbb{Z}^3}\frac{e^{i\ell\,k\cdot\eta}}{|x+\ell\,k|}=\frac{4\pi}{\ell^3}\sum_{k\in \mathbb{Z}^3}\frac{1}{\left|\frac{2\pi k}{\ell}-\eta\right|^2}\,e^{i\left(\frac{2\pi k}{\ell}-\eta\right)\cdot x}.
\end{equation}
It is $Q_\ell^*$-periodic w.r.t. $\eta$ and quasi-periodic with quasi-momentum $\eta$ w.r.t. $x$. 
For every $\gamma\in \Gamma_{\leq q}$, we now define the mean-field periodic Dirac operator 
\[
\D_\gamma=\fint^\oplus_{Q^*_\ell} \D_{\gamma,\xi}\,d\xi \quad \text{ with }\quad 
\D_{\gamma,\xi}:= \D_\xi-z\,G_\ell+\alpha V_{\gamma,\xi}\]
where
\begin{equation}\label{eq:def-V}
V_{\gamma,\xi}=\rho_{\gamma}\ast G_\ell-W_{\gamma,\xi}
\end{equation}
with
\begin{equation}\label{eq:convol-G}
\rho_{\gamma}\ast G_\ell(x)= \int_{Q_\ell}G_\ell(y-x)\,\rho_{\gamma}(y)\,dy= \widetilde{\Tr}_{L^2}[G_\ell(\cdot-x)\,\gamma]
\end{equation}
and
\begin{align}\label{eq:op-W}
    W_{\gamma,\xi}\psi_\xi(x)=\fint_{Q_\ell^*}d\xi'\int_{Q_\ell}W_\ell^\infty(\xi'-\xi,x-y)\,\gamma_{\xi'}(x,y)\,\psi_\xi(y)\,dy.
\end{align}
(In \eqref{eq:convol-G} we keep the notation $\cdot\ast \cdot$ for the convolution of periodic functions on $Q_\ell$, and we define 
\[
\widetilde{\Tr}_{L^2}(\gamma):=\fint_{Q_\ell^*}\Tr_{L^2_\xi}(\gamma_\xi)\,d\xi,
\]
where the $\;\widetilde\quad$ reminds us that $\gamma$ is not trace-class on $L^2(\mathbb{R}^3)$.) 

Then we can rewrite the DF functional as follows
\begin{align}\label{eq:DFfunctional}
   \mathcal{E}(\gamma)=\widetilde{\Tr}_{L^2} [\D_{\gamma}\gamma ]-\frac{\alpha}{2} \widetilde{\Tr}_{L^2}[V_{\gamma}\gamma]. 
\end{align}
In the standard DF theory, the system of units is such that  $m=c=\hbar=1$ and $z$ should be replaced by $\alpha z$, where $m$ is the mass of the electron, $c$ the speed of light, and $\hbar$ the Planck constant. Consequently, the fine structure constant $\alpha\approx \frac{1}{137}$. 

In this paper, we rather consider the non-relativistic regime and the weak electron-electron interaction regime. The \textit{non-relativistic regime} corresponds to the case when $q$ and $z$ are fixed and $c\gg1$, whereas the \textit{weak coupling regime} means  $q$ and $z$ fixed and  $\alpha\ll1$. In both cases, we assume in the following, without loss of generality, that $c\geq 1$ and $\alpha\leq 1$. For future convenience, we denote 
\begin{align*}
\alpha_c:=\frac{\alpha}{c}\qquad\textrm{and}\qquad z_c:=\frac{z}{c}.
\end{align*}

\subsection{The periodic DF ground-state  energy}
Let 
\[
P_{\gamma}^{\pm}=\fint^\oplus_{Q^*} P_{\gamma,\xi}^{\pm}\,d\xi \quad \text { with }\quad  P_{\gamma,\xi}^{\pm}:=\mathbbm{1}_{\mathbb{R}_{\pm}}(\D_{\gamma,\xi})
\]
denote the projection onto the positive and negative spectrum of $\D_\gamma$, respectively. Note that by definition $P_{0,\xi}^{\pm}=\mathbbm{1}_{\mathbb{R}_{\pm}}(\D_\xi- zG_\ell)$. We define 
\[
\Gamma_{ q}^+:=\Set*{\gamma\in \Gamma_{ q}\given \gamma=P^+_\gamma \gamma P^+_\gamma},\quad \Gamma_{\leq q}^+:=\Set*{\gamma\in \Gamma_{\leq q}\given \gamma=P^+_\gamma \gamma P^+_\gamma}
\]
and the ground-state energy 
\begin{equation}\label{eq:Iq}
    I_q:=\inf_{\gamma\in \Gamma_q}\mathcal{E}(\gamma).
\end{equation}
Existence of a ground-state has been proved in \cite{crystals} under the following. 
\begin{assumption}\label{ass:1}
Let $q\in \mathbb{N}^+$, $z\in \mathbb{R}^+$, and   $\kappa(\alpha,c):=(C_Gz_c+C_{EE} \alpha_c q)$\footnote{In \cite{crystals}, $q^+=\max\{q,1\}$ is introduced. As $q\in\mathbb{N}^+$, in this paper we have $q^+=q$.}. We assume that
\begin{enumerate}
    \item $\kappa(\alpha,c)<1-\frac{\alpha_c}{2} C_{EE}q$;
   \item There exists $\rho>0$ such that :
    \begin{align*}
        \begin{cases}
        1<\rho< \alpha_c^{-1}(1-\kappa(\alpha,c))^{1/2}\lambda_0^{1/2}(\alpha,c);\\
        \sqrt{\max\{(1-\kappa(\alpha,c)-\frac{\alpha_c}{2} C_{EE}q)^{-1}c^{-2}c^*(q+1)\,q,1\}q}<\rho,
    \end{cases}
    \end{align*}
\end{enumerate}
where $C_{EE}:=C_{EE}(\ell)$ is a constant given by Lemma \ref{lem:hardy}, $\lambda_0(\alpha,c)$ is given by Lemma \ref{lem:ope}, and $c^*(q+1)$ is defined in Lemma \ref{lem:prop-spec}.
\end{assumption}
\begin{theorem}[Existence of minimizers {\cite[Theorem 2.5]{crystals}}] \label{3.th:2.7'}
Under Assumption \ref{ass:1} above, the minimum problem $I_q$ admits a mini\-mi\-zer $\gamma_*\in\Gamma_{ q}^+$. Furthermore, $\fint_{Q_\ell^*} \Tr_{L^2_\xi}(\gamma_{{*},\xi})\,d\xi=q$ and $\gamma_*$ solves the  following nonlinear self-consistent equation:
\begin{equation}\label{3.eq:2.21}
\begin{aligned}
\gamma_*=\mathbbm{1}_{[0,\nu)}(\D_{\gamma_*})+\delta
\end{aligned}
\end{equation}
where $0\leq \delta \leq \mathbbm{1}_{\{\nu\}}(D_{\gamma_*})$ and  $\nu$ is the Lagrange multiplier due to the charge constraint $\fint_{Q_\ell^*} \Tr_{L^2_\xi}(\gamma_{\xi})\,d\xi= q$ satisfying  $0\leq \nu\leq c^*(q+1)$. 
\end{theorem}

\subsection{Main result}
We assume the following.
\begin{assumption}\label{ass:2}
Let $q\in \mathbb{N}^+$ and $z\in \mathbb{R}^+$ be fixed. We assume that $0<\alpha\leq 1$ and $c\geq 1$ are chosen such that
\begin{enumerate}
\item Assumption \ref{ass:1} is satisfied;
\item $c\geq 2(C_Gz+C_{EE}q)$;
\item $\alpha\leq\frac{4\pi}{C_{\rm cri}(q,z)}\,c$ where $C_{\rm cri}(q,z)$ is a large enough constant independent of $\alpha$ and $c$ that is given in Lemma \ref{lem:contra2}.
\end{enumerate}
\end{assumption}
In the non-relativistic limit $c\to+\infty$, we have $\kappa(\alpha,c)\to 0$. In the same manner, in the weak coupling limit $\alpha\to 0^+$, we have $\kappa(\alpha,c)\to C_G z_c$. Then it is easy to see that Assumption \ref{ass:2} is  satisfied in both situations whenever  $2(C_Gz+C_{EE}q)<c$. 

The purpose of this paper is to show that, under the above conditions on $\alpha$ and $c$, we have the following. 
\begin{theorem}[Properties of the last shell]\label{th:prop}
Under Assumption \ref{ass:2}, for any minimizer $\gamma_*$ of \eqref{eq:Iq}, the density matrix $\delta$ given in \eqref{3.eq:2.21} satisfies either $\delta=0$ or $\delta=\mathbbm{1}_{\{\nu\}}(\D_{\gamma_*})$. As a result, $\gamma_*$ is a projector.
\end{theorem}
\begin{remark}
 Assumption \ref{ass:1} is not empty. According to \cite[Remark 2.7]{crystals}, for $c=1$, $\alpha=\frac{1}{137}$ and $\ell\approx 1000$, Assumption \ref{ass:1} is satisfied for $q=z\leq 17$. Condition 2 in Assumption \ref{ass:2} is used to guarantee that $\kappa(\alpha,c)$ is away from $1$ uniformly. As we consider the case $c\gg 1$, this is not a problem even if we choose $c\geq 2(C_Gz+C_{EE}q)$. 
\end{remark}

\begin{remark}
We point out the fact that this is the first time that an explicit upper bound on $\frac{\alpha}{c}$ is given which ensures  that any DF minimizer is a projector. (In the DF model for atoms and molecules \cite{meng}, the same result is not available.)

 As in the HF model, these two cases rely on totally different argument: For atoms and molecules, the positive definiteness of the nonlinear term is used; while for crystals, the proof is based on a careful analysis of the singularity w.r.t. $\xi\in Q_\ell^*$ of the nonlinear term. Compared with the former, the singularity of the latter provides a quantitative estimate (i.e., \eqref{eq:4.10}) independent of $\alpha$ and $c$. This estimate gives the upper bound $C_{\rm cri}(q,z)$. 
\end{remark}

\begin{remark} 
 When $\alpha=0$, the nonlinear term $V_{\gamma}$ disappears. In this case, the proof that any minimizer is a projector relies on the absolute continuity of the spectrum of operators. In the non-relativistic case, the spectrum of the periodic operator $-\frac{1}{2}\Delta-zG$ is absolutely continuous \cite{thomas1973time}: One can use \cite[Theorem 1.9, Chapter 7.1]{kato2013perturbation} since $G_\ell(-\Delta)^{-1}$ is compact. On the contrary, the operator $G_\ell \,\D^{-1}$ is not compact, and we do not know whether the spectrum of the periodic Dirac-Coulomb operator $\mathcal{D}-zG$ is absolutely continuous. 
\end{remark}

The proof of  Theorem \ref{th:prop} relies on the ideas of Ghimenti and Lewin~\cite{ghimenti2009properties}. Their proof is based on the local convexity and the second-order expansion of the periodic HF functional on the constraint set. In the DF case, we are convinced that the constraint set $\Gamma^+_q$ is not convex and we are not able to prove that it is closed for the weak-$*$ topology. Both observations lead to the failure of the direct study of the second-order expansion of the periodic DF model on $\Gamma_q^+$. 

Following \cite{crystals}, instead of studying directly the DF model on $\Gamma_q^+$, we consider a penalized DF model (i.e., $\mathcal{E}(\gamma)-\epsilon_P\widetilde{\Tr}_{L^2}[\gamma]$) on $\Gamma_{\leq q}^+$. Then, using the retraction technique developed by S\'er\'e \cite{Ser09} (see also Catto-Meng-Paturel-S\'er\'e \cite{crystals} for crystals), we construct a retraction mapping $\theta$ onto a convex open subset $\mathcal{U}_R\subset \Gamma_{\leq q}^+$, such that any minimizer of $I_q$ is situated in $\mathcal{U}_R$, and such that we have $\theta(\gamma)\in \Gamma_{\leq q}^+$ for any $\gamma\in \mathcal{U}_R$.  As shown in \cite[Corollary 3.3]{crystals}, the two models are equivalent. Therefore,  we investigate in this paper the second-order expansion of the new  DF functional $\gamma \mapsto \mathcal{E}(\theta(\gamma))-\epsilon_P \widetilde{\Tr}_{L^2}[\theta(\gamma)]$ on $\mathcal{U}_R$. This property can be obtained by mimicking the proof of Meng \cite{meng}.

With these results in hand, Theorem~\ref{th:prop} is obtained  by following the lines of the proof in  \cite{ghimenti2009properties}.

\textit{Organisation of this paper:} In Section \ref{sec:3}, we first study the penalized DF model on $\Gamma_{\leq q}^+$, the retraction $\theta$, and the second-order expansion for the new DF functional. In Section \ref{sec:4}, we adapt the ideas of \cite{ghimenti2009properties} to prove Theorem \ref{th:prop}. In section \ref{sec:retra}, we recall the proof of Lemma \ref{lem:retra}. In Section \ref{sec:5}, we give the details about the proof of the second-order expansion for the new DF functional based on new functional space $\mathcal{Y}_c$. Some useful inequalities and  bounds on the eigenfunctions of DF operators are studied in Appendix.
\section{Relaxed DF energy and its second-order expansion}\label{sec:3}

It is shown in \cite{crystals} that, under Assumption \ref{ass:1},  existence of minimizers for problem $I_q$ is equivalent to the following penalized minimum problem for some $\epsilon_P$ large enough:
\begin{align}\label{eq:min1}
    I_{\leq q}:=\inf_{\gamma\in \Gamma_{\leq q}^+}(\mathcal{E}(\gamma)-\epsilon_{P} \widetilde{\Tr}_{L^2}[\gamma])+\epsilon_P q.
\end{align}
More precisely, we have the following. 
\begin{theorem}[Minimizers of the penalized problem {\cite[Theorem 3.2 and Corollary 3.3]{crystals}}]\label{th:penal}
We assume that Assumption \ref{ass:1} holds. Then there is a positive constant $\epsilon$ small enough independent of $\alpha,c$ such that for $\epsilon_P:=c^*(q+1)+\epsilon$ with $c^*(q+1)$ given in Lemma \ref{lem:prop-spec}, the penalized minimum problem $I_{\leq q}$ admits a minimizer $\gamma_*\in\Gamma_{\leq q}^+$ with  $\fint_{Q_\ell^*} \Tr_{L^2_\xi}(\gamma_{{*},\xi})\,d\xi=q$. In addition, $\gamma_*$ solves the following nonlinear self-consistent equation:
\begin{equation}\label{3.eq:2.21'}
\begin{aligned}
\gamma_*=\mathbbm{1}_{[0,\nu)}(\D_{\gamma_*})+\delta
\end{aligned}
\end{equation}
where $0\leq \delta \leq \mathbbm{1}_{\{\nu\}}(\D_{\gamma_*})$ and  $\nu$ is the Lagrange multiplier due to the charge constraint $\fint_{Q_\ell^*} \Tr_{L^2_\xi}(\gamma_{*,\xi})\,d\xi\leq q$ with  $0\leq \nu\leq c^*(q+1)$. 

Moreover, $I_q=I_{\leq q}$. As a result, any minimizer $\gamma_*$ of $I_{\leq q}$ is a minimizer of $I_q$ and vice versa. 
\end{theorem}
The proof of existence of minimizers of \eqref{eq:min1} relies on the construction of  a regular retraction $\theta$ defined  on an open subset $\mathcal{U}_R$ of $\Gamma_{\leq q}$ such that
\[
\theta:\quad \mathcal{U}_R \to \mathcal{U}_R\bigcap \Gamma_{\leq q}^+,\qquad \textrm{and}\qquad\theta(\mathcal{U}_R\bigcap \Gamma_{\leq q}^+)=\mathcal{U}_R\bigcap \Gamma_{\leq q}^+.
\] 
For $\gamma \in \Gamma_{\leq q}$, the retraction is defined by 
\begin{align}\label{eq:theta}
    \theta(\gamma):=\lim_{n\rightarrow +\infty} T^n(\gamma),
\end{align}
with
\[
T(\gamma)=P^+_{\gamma}\gamma P^+_{\gamma},\quad T^n(\gamma)=T(T^{n-1}(\gamma)),\quad T^0(\gamma)=\gamma.
\]
According to Floquet unitary transform, we have
    \begin{align*}
        T(\gamma)=\fint_{Q_\ell^*}^\oplus T_\xi(\gamma)\,d\xi=\fint_{Q_\ell^*}^\oplus P^+_{\gamma,\xi}\gamma_\xi P^+_{\gamma,\xi}\,d\xi.
    \end{align*}
With $\theta$ in hand, if a minimizer $\gamma_*$ of $I_{\leq q}$ is in $\mathcal{U}_R$, the penalized minimization problem $I_{\leq q}$ reduces to a simpler one :
\begin{align}\label{min-retra}
    I_{\leq q}=\min_{\gamma\in \mathcal{U}_{R}}E(\gamma)+\epsilon_P\,q,
\end{align}
where the new DF functional $E(\cdot)$ is defined by 
\begin{align}\label{eq:E}
    E(\gamma):=\mathcal{E}(\theta(\gamma))-\epsilon_P\widetilde{\Tr}_{L^2}[\theta(\gamma)].
\end{align}
The following lemma guarantees the existence of the retraction $\theta$ in our setting where the norm $Y$ in \cite[Proposition 5.8]{crystals} is replaced by $\mathcal{Y}_c$. 
\begin{lemma}\label{lem:retra}
Let $\kappa(\alpha,c)<1$ and $A(\alpha,c):=\frac{\alpha_c}{2}C_{EE}(1-\kappa(\alpha,c))^{-1}$. Given $1<R<\frac{1}{2A(\alpha,c)}$, let $M(\alpha,c):=\max\left(\frac{1+A(\alpha,c)q}{2},\frac{1}{1-2A(\alpha,c)R}\right)$, and let
\[
\mathcal{U}_R:=\Set*{\gamma\in \Gamma_{\leq q}\given \frac{1}{c}\max\{\|\gamma|\D|^{1/2}\|_{\mathfrak{S}_{1,1}},\|\gamma\|_{\mathcal{Y}}\}+\frac{M(\alpha,c)}{c^2}\|T(\gamma)-\gamma\|_{X_c\bigcap \mathcal{Y}_c}< R}.
\]
Then, $\mathcal{U}_R$ is a open subset in $\Gamma_{\leq q}$, $T$ maps continuously $\mathcal{U}_R$ into $\mathcal{U}_R$, and for any $\gamma\in\mathcal{U}_R$ the sequence $(T^n(\gamma))_{n\geq 0}$ converges to a limit $\theta(\gamma)\in\Gamma_q^+$. Moreover, for all $\gamma\in \mathcal{U}_R$,
\begin{align}
\|T^{n+1}(\gamma)-T^n(\gamma)\|_{X_c\bigcap \mathcal{Y}_c}\leq L(\alpha,c)\|T^{n}(\gamma)-T^{n-1}(\gamma)\|_{X_c\bigcap \mathcal{Y}_c}\label{eq:3.3}
\end{align}
and
\begin{align}
   \|\theta(\gamma)-T^n(\gamma)\|_{X_c\bigcap \mathcal{Y}_c}\leq \frac{L^n(\alpha,c)}{1-L(\alpha,c)}\|T(\gamma)-\gamma\|_{X_c\bigcap \mathcal{Y}_c},\label{eq:3.4}
\end{align}
with $0< L(\alpha,c)=2A(\alpha,c)R<1$.
\end{lemma}
This result is collected in \cite[Proposition 5.4]{crystals} and can be proved by adapting \cite[Proposition 5.8]{crystals} and \cite[Proposition 2.1]{Ser09} to our new functional framework $\mathcal{Y}$. For the reader's convenience, the proof of Lemma \ref{lem:retra} in this context is provided in Section \ref{sec:retra}.

The main ingredient of the proof of Theorem \ref{th:prop} is the following.
\begin{theorem}[Second-order expansion for the new DF energy] \label{approx}
Let $R>1$ be fixed, and let $\kappa(\alpha,c)<1$ and $L(\alpha,c)<1$ be given as in Lemma \ref{lem:retra}. Given any $\gamma\in \mathcal{U}_R\bigcap \Gamma_{\leq q}^+$ and $h\in X\bigcap Y$ such that $P^+_{\gamma}h P^+_{\gamma}=h$, and given any $t$ such that $\gamma+th\in \mathcal{U}_R$, we have
\begin{equation}\label{eq:second}
    \begin{aligned}
\MoveEqLeft    E(\gamma+th)=E(\gamma)+t\,\widetilde{\Tr}_{L^2} [(\D_{\gamma}-\epsilon_P)h]+\frac{\alpha t^2}{2}\widetilde{\Tr}_{L^2}[V_h h]+t^2\alpha_c^2  Err(t,\gamma,h)
\end{aligned}
\end{equation}
where $|Err(t,\gamma,h)|\leq (2+\frac{\epsilon_P}{c^2})\mathcal{N}_{\gamma+th}(h)$ with
\begin{align}\label{eq:error}
\mathcal{N}_\gamma(h)&=\frac{ C_{EE}^2}{2(1-\kappa(\alpha,c))^{2}\lambda_0(\alpha,c)}\fint_{Q_{\ell}^*}\|h\|_{X\bigcap \mathcal{Y}(\xi)}^2 \|\gamma_\xi\|_{\mathfrak{S}_1(\xi)}d\xi\notag\\
 &+(q\alpha_c^2+R\alpha_c^2 +\alpha_c) \frac{10C_{EE}^4}{(1-\kappa(\alpha,c))^4\lambda_0(\alpha,c)^{5/2}(1-L(\alpha,c))^2} \notag\\
  &\qquad\times\left(\frac{1}{c}\fint_{Q_{\ell}^*}\|h\|_{X\bigcap \mathcal{Y}(\xi)}\|\gamma_\xi|\D_\xi|^{1/2}\|_{\mathfrak{S}_1(\xi)}d\xi
+\frac{1}{c}\sup_{\xi\in Q_\ell^*}\fint_{Q_{\ell}^*}\frac{\|h\|_{X\bigcap \mathcal{Y}(\xi')}}{|\xi-\xi'|^2}d\xi'\right)^2.
\end{align}
\end{theorem}
The proof is provided in Section \ref{sec:5}.
\begin{remark}
For $q\in \mathbb{N}^+$, it is not difficult to see that
\[
\mathcal{N}_\gamma(h)\leq (R+q)\|h\|_{X\bigcap Y}
\]
by using $\gamma\in\mathcal{U}_R$ as in \cite{meng}. As a result, $Err(t,\gamma,h)$ is uniformly bounded with respect to $t$ if $\gamma+th\in \mathcal{U}_R$. However, in this paper, we can not simplify \eqref{eq:error}. We have to use this complicated formula to control the delicate behavior of $\alpha_c$ in Lemma  \ref{lem:contra2}. 
\end{remark}
We now embark on the proof of Theorem \ref{th:prop}.
\section{Properties of the last shell}\label{sec:4}
Now we are going to mimic the proof in \cite{ghimenti2009properties} to prove Theorem \ref{th:prop} with the help of Theorem \ref{approx}.

\subsection{Continuity of the eigenfunctions of \texorpdfstring{$\mathcal{D}_{\gamma_*,\xi}$}{}  w.r.t. \texorpdfstring{$\xi$}{}}
We recall that $W_{\ell}^\infty$ is defined by \eqref{3.eq:2.8}. We may separate the singularities of $W_\ell^\infty$ w.r.t. $\eta\in 2Q_\ell^*$ and $x\in 2Q_\ell$ as follows
\begin{equation}\label{3.eq:4.7'}
    W_\ell^\infty(\eta,x)=W_{\geq 2,\ell}(\eta,x)+W_{<2,\ell}(\eta,x),
\end{equation}
with 
\[
W_{\geq 2,\ell}(\eta,x)=\frac{4\pi}{\ell^3}\sum_{\substack{|k|_\infty\geq 2\\ k\in \mathbb{Z}^3}}\frac{1}{\left|\frac{2\pi k}{\ell}-\eta\right|^2}\,e^{i\left(\frac{2\pi k}{\ell}-\eta\right)\cdot x}
\]
and
\[
W_{<2,\ell}(\eta,x)=\frac{4\pi}{\ell^3}\sum_{\substack{|k|_\infty< 2\\ k\in \mathbb{Z}^3}}\frac{1}{\left|\frac{2\pi k}{\ell}-\eta\right|^2}\,e^{i\left(\frac{2\pi k}{\ell}-\eta\right)\cdot x},
\]
where $|k|_\infty:=\max\{|k_1|,|k_2|,|k_3|\}$.
It is convenient to study $\D_{\gamma_*}$ on the fixed Hilbert space $L^2_{\rm per}(Q_\ell)$ with $\gamma_*$ being a minimizer of \eqref{eq:min1}. Thus we introduce the unitary operator defined in each Bloch fiber by
\[
U_\xi:L^2_{\rm per}(Q_\ell)\to L^2_\xi(Q_\ell),\quad u\mapsto e^{i\xi\cdot x}u(\cdot).
\]
When $A_\xi$ is an operator on $L^2_\xi(Q_\ell)$, we  shall use the notation  $\widetilde{A}_\xi:=U_\xi^* A_\xi U_\xi$ for the corresponding operator on $L^2_{\rm per}(Q_\ell)$. With the operator $W_{\gamma_*,\xi}$ being defined by \eqref{eq:op-W}, we have the operator $\widetilde{W}_{\gamma_*,\xi}$ defined by its kernel
\begin{align*}
  \widetilde{W}_{\gamma_*,\xi}(x,y)&=[U_\xi^* W_{\gamma_*,\xi}U_\xi](x,y)\\
   &=\fint_{Q_\ell^*}\Bigg[\widetilde{W}_{\geq 2,\ell}(\xi-\xi',x-y)+\widetilde{W}_{<2,\ell}(\xi-\xi',x-y)\Bigg]\widetilde{\gamma}_{*,\xi'}(x,y)d\xi',
\end{align*}
where $\widetilde{W}_{\geq 2,\ell}(\eta,x)=e^{i\eta\cdot x} W_{\geq 2,\ell}(\eta,x)$, $\widetilde{W}_{<2,\ell}(\eta,x)=e^{i\eta\cdot x} W_{<2,\ell}(\eta,x)$ and $\widetilde{\gamma}_{*,\xi'}=e^{-i\xi'\cdot x}\gamma_{*,\xi'}e^{i\xi'\cdot y}$.

\begin{lemma}\label{lem:Holder}\footnotemark
Let $0< p< 1$. 
The family $(\widetilde{W}_{\gamma_*,\xi})_{\xi\in Q_\ell^*}$ is bounded and  $p$-H\"older continuous in $\mathcal{B}(L^2_{\rm per}(Q_\ell))$ w.r.t. $\xi\in Q_\ell^*$. Moreover, for every $\xi_1$ and $\xi_2$ in $Q_\ell^*$,
\begin{align*}
     \|\widetilde{W}_{\gamma_*,\xi_1}-\widetilde{W}_{\gamma_*,\xi_2}\|_{\mathcal{B}(L^2_{\rm per}(Q_\ell))}\leq C_p\,|\xi_1-\xi_2|^p,
\end{align*}
with $C_p$ being a positive constant independent of $\xi_1$ and $\xi_2$.
\end{lemma}
\footnotetext{It is claimed in \cite[Page 745]{catto2001thermodynamic} and \cite[Lemma 1]{ghimenti2009properties} that the function $f(\eta,x):=W_\ell^\infty(\eta,x)-e^{-i\eta\cdot x}G_\ell(x)-\frac{4\pi}{\ell^3}\frac{e^{-i\eta\cdot x}}{|\eta|^2}$ is harmonic which is not true.}
\begin{proof}
The boundedness of $\widetilde{W}_{\gamma_*.\xi}$ follows from \eqref{eq:W-Y}. Indeed, 
\begin{align*}
   \supess_{\xi\in Q_{\ell}^*}\|\widetilde{W}_{\gamma_*,\xi}\|_{\mathcal{B}(L^2_{\rm per}(Q_\ell))}&=\supess_{\xi\in Q_{\ell}^*}\|U^*_\xi W_{\gamma_*,\xi}U_\xi\|_{\mathcal{B}(L^2_{\rm per}(Q_\ell))}=\|W_{\gamma_*}\|_Y\leq C\|\gamma_*\|_{X\bigcap \mathcal{Y}}.
\end{align*}
Using the fact that $\widetilde{W}_{\geq 2,\ell}(\eta,x)=e^{i\eta\cdot x} W_{\geq 2,\ell}(\eta,x)$, we have
\begin{align*}
\MoveEqLeft    \widetilde{W}_{\geq 2,\ell}(\xi_1-\xi',x-y)-\widetilde{W}_{\geq 2,\ell}(\xi_2-\xi',x-y)\\
&=\int_{0}^1  i  (\xi_1-\xi_2)\cdot(x-y)e^{i(\xi_2+t(\xi_1-\xi_2)-\xi')\cdot (x-y)} W_{\geq 2,\ell}(\xi_2+t(\xi_1-\xi_2)-\xi',x-y)\\
&\qquad\qquad+ [e^{i(\xi_2+t(\xi_1-\xi_2)-\xi')\cdot (x-y)}\nabla_\xi W_{\geq 2,\ell}(\xi_2+t(\xi_1-\xi_2)-\xi',x-y)]\cdot (\xi_1-\xi_2)d t.
\end{align*}
For any $u,v\in L^2_{\rm per}(Q_\ell)$, using \eqref{3.eq:4.7'} we have
\begin{align*}
    \MoveEqLeft \left|\left(u,(\widetilde{W}_{\gamma_*,\xi_1}-\widetilde{W}_{\gamma_*,\xi_2})v\right)_{L^2_{\rm per}}\right|\\
    &\leq C|\xi_1-\xi_2|\sup_{\xi\in 2Q_\ell^*}\int_{Q_\ell}\int_{Q_\ell}\fint_{Q_\ell^*} |\widetilde{W}_{\geq 2,\ell}(\xi-\xi',x-y)||\widetilde{\gamma}_{*,\xi'}(x,y)||u(y)||v(x)|d\xi'dxdy\\
    &\quad+ |\xi_1-\xi_2|\fint_{Q_\ell^*}d\xi'\sup_{\xi\in 2Q_\ell^*}\|\nabla_{\xi}W_{\geq 2,\ell}(\xi-\xi')\widetilde{\gamma}_{*,\xi'}\|_{L^2(Q_\ell\times Q_\ell)}\|u\|_{L^2_{\rm per}(Q_\ell)}\|v\|_{L^2_{\rm per}(Q_\ell)}\\
    &\quad +\frac{4\pi}{\ell^3}\sum_{\substack{k\in \mathbb{Z}^3\\|k|_\infty<2}}\fint\limits_{Q^*_\ell}d\xi'\left|\frac{1}{\big|\xi'-\xi_1-\frac{2\pi k}{\ell}\big|^2}-\frac{1}{\big|\xi'-\xi_2-\frac{2\pi k}{\ell}\big|^2}\right|\left|\left(\widetilde{\gamma}_{*,\xi'}U_{2\pi k}u, U_{2\pi k}v\right)_{L^2_{\rm per}}\right| .
\end{align*}
It is shown in \cite[Lemma B.2 and Lemma B.3]{crystals} that
\begin{align*}
    \sup_{\xi\in2Q_\ell^*}\||W_{\geq 2,\ell}|^{1/2}(\xi,\cdot)|\D_\xi|^{-1/2}\|_{\mathcal{B}(L^2_{\xi})}\leq C ,\quad\sup_{\xi\in2Q_\ell^*}\|\nabla_\xi W_{\geq 2,\ell}(\xi,\cdot)\|_{L^\infty(Q_\ell)}\leq C.
\end{align*}
Then as $\|\gamma\|_Y\leq 1$, from \eqref{eq:W-Y} and arguing as for \cite[Eq. (B.9)]{crystals}, we get
\begin{align*}
    \MoveEqLeft \left|\left((\widetilde{W}_{\gamma_*,\xi_1}-\widetilde{W}_{\gamma_*,\xi_2})u,v\right)_{L^2_{\rm per}}\right|\leq C|\xi_1-\xi_2|\|\gamma_*\|_{X}\\
    &\quad+C|\xi_1-\xi_2|\left(\fint_{Q_\ell}^*\int_{Q_\ell\times Q_\ell}|\widetilde{\gamma}_{*,\xi'}(x,y)|^2dxdyd\xi'\right)^{1/2}\|u\|_{L^2_{\rm per}(Q_\ell)}\|v\|_{L^2_{\rm per}(Q_\ell)}\\
    &\quad+C\sum_{\substack{k\in \mathbb{Z}^3\\|k|_\infty<2}}\fint\limits_{Q^*_\ell}d\xi'\left|\frac{1}{\big|\xi'-\xi_1-\frac{2\pi k}{\ell}\big|^2}-\frac{1}{\big|\xi'-\xi_2-\frac{2\pi k}{\ell}\big|^2}\right|\|u\|_{L^2_{\rm per}(Q_\ell)}\|v\|_{L^2_{\rm per}(Q_\ell)}.
\end{align*}
Observe that
\[
\fint_{Q_\ell^*}\int_{Q_\ell\times Q_\ell}|\widetilde{\gamma}_{*,\xi'}(x,y)|^2dxdyd\xi'=\fint_{Q_\ell^*}\Tr_{L^2_{\rm per}(Q_\ell)}\widetilde{\gamma}_{*,\xi'}^2d\xi'\leq \fint_{Q_\ell^*}\Tr_{L^2_{\rm per}(Q_\ell)}\widetilde{\gamma}_{*,\xi'}\,d\xi'=q.
\]
Then
\[
\|\widetilde{W}_{\gamma,\xi_1}-\widetilde{W}_{\gamma,\xi_2}\|_{\mathcal{B}(L^2_{\rm per}(Q_\ell))}\leq C|\xi_1-\xi_2|+C\sum_{\substack{k\in \mathbb{Z}^3\\|k|_\infty<2}}\fint\limits_{Q^*_\ell}d\xi'\left|\frac{1}{\big|\xi'-\xi_1-\frac{2\pi k}{\ell}\big|^2}-\frac{1}{\big|\xi'-\xi_2-\frac{2\pi k}{\ell}\big|^2}\right|.
\]
The H\"older continuity follows since the function $\xi \mapsto \fint_{Q_\ell^*}\frac{1}{|\xi-\xi'|^2}d\xi'$ is $p$-H\"older continuous w.r.t. $\xi$ for any $0<p<1$.
\end{proof}
We also need the following, which adapts \cite[Lemma 4]{ghimenti2009properties} to the periodic Dirac operator.
\begin{lemma}
Let $\kappa(\alpha,c)<1$. Consider $\Omega$ an open subset of $Q_\ell^*$. For any $\epsilon>0$, let $K$ be a compact set in $\mathbb{C}$ such that $\inf_{\xi\in \Omega} d(K;\sigma(\widetilde{\D}_{\gamma_*,\xi}))\geq \epsilon>0$. Then we have
\begin{enumerate}
    \item $\D(\widetilde{\D}_{\gamma_*,\xi}-z)^{-1}$ is bounded on $L^2_{\rm per}(Q_\ell)$, uniformly w.r.t. $\xi\in \Omega$ and $z\in K$;
    \item The map $\xi\mapsto \D(\widetilde{\D}_{\gamma_*,\xi}-z)^{-1}$ is H\"older continuous w.r.t. $\xi\in \Omega$  with values in $\mathcal{B}(L^2_{\rm per}(Q_\ell))$, uniformly in $z\in K$.
\end{enumerate}
\end{lemma}
\begin{proof}
As $\widetilde{\D}_\xi=\D +c\sum_{k=1}^3\alpha_k\xi_k$ and $\xi\in \Omega\subset Q_\ell^*$, we have
\begin{align*}
\MoveEqLeft    \|\D(\widetilde{\D}_{\gamma_*,\xi}-z)^{-1}\|_{\mathcal{B}(L^2_{\rm per}(Q_\ell))}\\
&\leq \|\widetilde{\D}_\xi(\widetilde{\D}_{\gamma_*,\xi}-z)^{-1}\|_{\mathcal{B}(L^2_{\rm per}(Q_\ell))}+c\||\xi| (\widetilde{\D}_{\gamma_*,\xi}-z)^{-1}\|_{\mathcal{B}(L^2_{\rm per}(Q_\ell))}\\
    &\leq\|\D_\xi(\D_{\gamma_*,\xi}-z)^{-1}\|_{\mathcal{B}(L^2_{\xi}(Q_\ell))}+\frac{C}{\epsilon}.
\end{align*}
Thanks to Lemma \ref{lem:ope}, we have,  for all $\xi\in \Omega$ and $z\in K$,
\begin{align*} \|\D_\xi(\D_{\gamma_*,\xi}-z)^{-1}\|_{\mathcal{B}(L^2_{\xi}(Q_\ell))}&\leq \frac{1}{1-\kappa(\alpha,c)}\|\D_{\gamma_*,\xi}(\D_{\gamma_*,\xi}-z)^{-1}\|_{\mathcal{B}(L^2_{\xi}(Q_\ell))}\\
    &\leq \frac{1}{1-\kappa(\alpha,c)}+\frac{|z|}{
    \epsilon(1-\kappa(\alpha,c))},
\end{align*}
which implies
\[
\|\D(\widetilde{\D}_{\gamma_*,\xi}-z)^{-1}\|_{\mathcal{B}(L^2_{\rm per}(Q_\ell))}\leq C,\quad \textrm{uniformly w.r.t. }\xi\in \Omega,\,z\in K.
\]
The H\"older continuity follows from the fact
\begin{align*}
\MoveEqLeft    \D[(\D_{\gamma_*,\xi}-z)^{-1}-(\D_{\gamma_*,\xi'}-z)^{-1}]\\
&= \D(\D_{\gamma_*,\xi}-z)^{-1}[\D_{\gamma_*,\xi}-\D_{\gamma_*,\xi'}](\D_{\gamma_*,\xi'}-z)^{-1}\\
&=[\D(\D_{\gamma_*,\xi}-z)^{-1}]\left[-ic\sum_{k=1}^3\alpha_k(\xi_k-\xi_k')-\widetilde{W}_{\gamma_*,\xi}+\widetilde{W}_{\gamma_*,\xi'}\right](\D_{\gamma_*,\xi'}-z)^{-1}.
\end{align*}
Here $\D(\D_{\gamma_*,\xi}-z)^{-1}$ and $(\D_{\gamma_*,\xi'}-z)^{-1}$ are bounded on $L^2_{\rm per}(Q_\ell)$, and $\xi\mapsto \widetilde{W}_{\gamma_*,\xi}$ is H\"older continuous as shown in Lemma \ref{lem:Holder}.
\end{proof}
We denote by $\lambda_k(\xi)$ the positive eigenvalues of $\widetilde{\D}_{\gamma_*,\xi}$ counted with multiplicity, which are the same as that the ones of $(\D_{\gamma_*,\xi})_\xi$ since $U_\xi$ is unitary. We may assume that the $\lambda_k(\xi)$'s are counted in nondecreasing order: $\lambda_1(\xi)\leq \lambda_2(\xi)\leq \cdots$.
 
Arguing as in \cite[Lemma 5 and Lemma 6]{ghimenti2009properties} and using Cauchy's formula \cite{kato2013perturbation}, we can immediately obtain the following lemma. 
\begin{lemma} Let $\kappa(\alpha,c)<1$. For any $k\in \mathbb{N}$,
the positive eigenvalues $\lambda_k(\xi)$ of $\D_{\gamma_*,\xi}$ are H\"older continuous  with respect to $\xi\in Q_\ell^*$. Furthermore, there exists $R>0$ independent of $\alpha$ and $c$ such that let $\Omega$ be an open subset of $Q_\ell^*$ and $I=(a,b)$ an interval of $\mathbb{R}^+$ such that $\sigma(\widetilde{\D}_{\gamma_*,\xi})\bigcap \{a,b\}=\emptyset$ for all $\xi\in \Omega$. Then, the map
\[
\xi\in\Omega\mapsto \D\mathbbm{1}_{I}(\widetilde{\D}_{\gamma_*,\xi})\in \mathcal{B}(L^2_{\rm per}(Q_\ell))
\]
is H\"older continuous. In particular, we can find an orthonormal basis $(u_1(\xi),\cdots,u_K(\xi))$ of the range of $\mathbbm{1}_I(\widetilde{\D}_{\gamma,\xi})$ such that $\xi\in \Omega\mapsto u_k(\xi)\in H^1_{\rm per}(Q_\ell)$ is H\"older continuous  w.r.t. $\xi\in \Omega$.
\end{lemma}

\subsection{The Fermi level is either empty or totally filled.}
Now we argue by contradiction and assume as in \cite{ghimenti2009properties} that $\delta\neq 0$ and $\delta\neq \mathbbm{1}_{\{\nu\}}(\D_{\gamma_{*}})$ and that the Lagrange multiplier $\nu$ defined in Theorem~\ref{3.th:2.7'} is an eigenvalue of $\D_{\gamma_*}$ (otherwise, $\delta_\xi=0$ for almost every $\xi\in Q_\ell^*$). Then,  
\[
|\{\xi\in Q_\ell^*; \exists k\geq 1, \lambda_k(\xi)=\nu\}|\neq 0.
\]
As in \cite{ghimenti2009properties}, we have the following lemma which is obtained without difficulty by replacing the HF operator $H_{\gamma_{*}}$ by the DF operator $\D_{\gamma_*}$ in the proof of \cite[Lemma 7]{ghimenti2009properties}.  

\begin{lemma}\label{lem:7.6}
Let $\kappa(\alpha,c)<1$. Assume that $\nu>0$ is an eigenvalue of $\D_{\gamma_*}$ such that  $\delta\neq 0$ and $\delta\neq \mathbbm{1}_{\{\nu\}}(\D_{\gamma_*})$. Then, there exists a constant $\epsilon>0$, a Borel set $\omega\subset Q_\ell^*$ with $|\omega|\neq 0$ and two continuous functions $\xi\in\omega\mapsto u(\xi)\in H^1_{\rm per}(Q_\ell)$ and $\xi\in \omega\mapsto u'(\xi)\in H^1_{\rm per}(Q_\ell)$ such that
\begin{align*}
   u(\xi),u'(\xi)\in \ker((\widetilde{\D}_{\gamma_*,\xi})_\xi-\nu),\quad \|u(\xi)\|_{L^2_{\rm per}(Q_\ell)}=\|u'(\xi)\|_{L^2_{\rm per}(Q_\ell)}=1, 
\end{align*}
and, denoting $\psi(\xi)=U^*_\xi u(\xi)$ and $\psi'(\xi)=U^*_\xi u'(\xi)$,
\begin{align}
    0\leq \delta_\xi+t\left|\psi(\xi)\right>\left<\psi(\xi)\right|-t'\left|\psi'(\xi)\right>\left<\psi'(\xi)\right|\leq 1
\end{align}
on $L^2_\xi(Q_\ell)$ for all $\xi\in \omega$ and all $t,t'\in [0,\epsilon)$, where $\left|\psi\right>\left<\psi\right|$ denotes the projector on the vector space spanned by the function $\psi$.
\end{lemma}
\begin{remark}
The hypothesis of $q\in \mathbb{N}^+$ is required in the proof of this lemma.
\end{remark}
The desired  contradiction is based on the positivity of the second-order term shown in Theorem \ref{approx}.
\begin{lemma}\label{lem:contra0}
Let $0<\alpha\leq 1$ and $c\geq 1$. Assume that $\nu$ is an eigenvalue of $\D_{\gamma_*}$ given by \eqref{3.eq:2.21'}, that $\delta\neq 0$ and $\delta\neq \mathbbm{1}_{\{\nu\}}(\D_{\gamma_*})$, and let $\omega$, $\psi(\xi)$ and $\psi'(\xi)$ be as in Lemma \ref{lem:7.6}. Then, there exists $R>1$ such that if $\alpha$ and $c$ satisfy
    \begin{align}\label{eq:alpha-c-R}
  c\geq 2(C_Gz+C_{EE}q),\qquad \alpha_c\leq \frac{1}  {4C_{EE}R},
\end{align}
we have 
    \begin{align*}
        \gamma_*, \gamma_*+th\in \mathcal{U}_R,\quad \textrm{ provided }\quad |t|\leq \min\{\frac{\epsilon}{\max\{\|\eta\|_{L^\infty},\|\eta'\|_{L^\infty}\}},1\}
    \end{align*}
where the periodic density matrix $h=\fint_{Q_\ell^*}^\oplus h_\xi d\xi$ is  defined by
\[
h_\xi=\eta(\xi)\left|\psi(\xi)\right>\left<\psi(\xi)\right|-\eta'(\xi)\left|\psi'(\xi)\right>\left<\psi'(\xi)\right|,
\]
with $\eta$ and $\eta'\in L^\infty(\omega,\mathbb{R}^+)$ satisfying $\int_{\omega}\eta=\int_{\omega}\eta'$ and $\psi|_{Q_\ell^*\setminus \omega}=\psi'|_{Q_\ell^*\setminus \omega}=0$.

Furthermore, we also have
\begin{align}\label{eq:4.9}
    \frac{\alpha }{2}\widetilde{\Tr}_{L^2}[V_h h]+\alpha_c^2  Err(t,\gamma_*,h)\geq 0.
\end{align}
Here $Err(t,\gamma_*,h)$ is defined in Theorem \ref{approx} and satisfies
\begin{align*}
    |Err(t,\gamma_*,h)|\leq (2+\frac{\epsilon_P}{c^2})\mathcal{N}_{\gamma_*+th}(h).
\end{align*}
\end{lemma}
\begin{proof}
In this proof, we denote by $C$ several different constants independent of $\alpha$ and $c$.

Let $t_0=\min\{\frac{\epsilon}{\max\{\|\eta\|_{L^\infty},\|\eta'\|_{L^\infty}\}},1\}$ where $\epsilon$ is defined in Lemma \ref{lem:7.6}. Let 
\begin{align}\label{eq:Rcontra}
    R=(K+C_M)q+C_{\mathcal{Y}}
\end{align}
with $K,\; C_{\mathcal{Y}}$ given by Lemma \ref{lem:unibound} and \eqref{eq:T-Y} respectively and $C_M$ given below by \eqref{eq:CM}. Before going further, we shall point out that \eqref{eq:alpha-c-R} implies that
\begin{align*}
    \kappa(\alpha,c)\leq \frac{C_Gz+C_{EE}\alpha q}{c} \leq \frac{1}{2}
\end{align*}
and
\begin{align*}
    L(\alpha,c)=2A(\alpha,c)R=\frac{\alpha_c C_{EE}}{(1-\kappa(\alpha,c))}R\leq \frac{1}{2}.
\end{align*}

We first verify that $\gamma_*\in \mathcal{U}_R$. By Lemma \ref{lem:prop-spec} and Lemma \ref{lem:unibound}, as $\kappa(\alpha,c)\leq \frac{1}{2}$,
\begin{equation}\label{eq:4.10'}
    \frac{1}{c}\|\gamma_*|\D|^{1/2}\|_{\mathfrak{S}_{1,1}}\leq \frac{1}{c}\|\gamma_*\|_{\mathfrak{S}_{1,1}}^{1/2}\|\gamma_*\|_{X_c}^{1/2}\leq\|\gamma_*\|_{\mathfrak{S}_{1,1}}^{1/2}\|\gamma_*\|_{X}^{1/2}\leq Kq
\end{equation}
and
\[
\|\gamma_*\|_{\mathcal{Y}}\leq \|\gamma_*\|_{Y}\fint_{Q_\ell^*}\frac{1}{|\xi|^2}d\xi\leq C_{\mathcal{Y}}.
\]
Thus,
\begin{align*}
    \MoveEqLeft \frac{1}{c}\max\{\|\gamma_*|\D|^{1/2}\|_{\mathfrak{S}_{1,1}},\|\gamma_*\|_{\mathcal{Y}}\}\!+\!\frac{M}{c^2}\|T(\gamma_*)\!-\!\gamma_*\|_{X_c\bigcap \mathcal{Y}_c}\!\\
    &\qquad\qquad\qquad=\frac{1}{c}\max\{\|\gamma_*|\D|^{1/2}\|_{\mathfrak{S}_{1,1}},\|\gamma_*\|_{\mathcal{Y}}\}\leq\! Kq+C_{\mathcal{Y}}\leq R.
\end{align*}
This implies $\gamma_*\in \mathcal{U}_R$.

Now we turn to study $\gamma_*+th\in \mathcal{U}_R$. For any $t\in [-t_0,t_0]$,  $\gamma_*+th\in \Gamma_q$ and $0\leq \gamma_*+th\leq \mathbbm{1}_{(0,c^2+\Sigma(q+1)]}(\D_{\gamma_*})$. Repeating the above estimate yields
\[
\frac{1}{c}\|(\gamma_*+th)|\D|^{1/2}\|_{\mathfrak{S}_{1,1}}\leq Kq,\quad  \|\gamma_*+th\|_{\mathcal{Y}}\leq C_{\mathcal{Y}}.
\]
Note that by the definition of $h$, we have $P^+_{\gamma_*}(\gamma_*+th)P^+_{\gamma_*}=\gamma_*+th$. Then we infer from Lemma \ref{lem:T-I} below that
\begin{align*}
    \frac{1}{c^2}\|T(\gamma_*+th)-\gamma_*-th\|_{X_c}&\leq Ct\frac{\alpha_c}{c^2} \fint_{Q_{\ell}^*}\|h\|_{X\bigcap \mathcal{Y}(\xi)}\|\gamma_{*,\xi}|\D_\xi|^{1/2}\|_{\mathfrak{S}_1(\xi)}d\xi\\
    &\leq C\frac{\alpha}{c^3}t\|h\|_{X\bigcap Y}\|\gamma_*|\D|^{1/2}\|_{\mathfrak{S}_{1,1}}\leq Cq\frac{\alpha}{c^2}
\end{align*}
and
\begin{align*}
    \frac{1}{c^2}\|T(\gamma_*+th)-\gamma_*-th\|_{\mathcal{Y}_c}&\leq Ct\frac{\alpha_c}{c^2} \|h\|_{X\bigcap \mathcal{Y}}\leq C\frac{\alpha}{c^2}
\end{align*}
since according to Lemma \ref{lem:unibound} we have
\[
t\|h\|_{X\bigcap \mathcal{Y}(\xi)}\leq Ct^2\|h\|_{X\bigcap Y}\leq C t_0\max\{\|\eta\|_{L^\infty},\|\eta'\|_{L^\infty}\}\leq C.
\]
By definition of $M(\alpha,c)$ in Lemma \ref{lem:retra}, it is easy to see that $M(\alpha,c)\leq C$. Thus there exists a constant $C_{M}>0$ independent of $\alpha$ and $c$ such that
\begin{align}\label{eq:CM}
    \frac{M(\alpha,c)}{c^2}\|T(\gamma_*+th)-\gamma_*-th\|_{X_c\bigcap \mathcal{Y}_c}\leq C_{M}q\frac{\alpha}{c^2}.
\end{align}
Now we conclude that 
\begin{align}\label{eq:R}
    \MoveEqLeft \frac{1}{c}\max\{\|(\gamma_*+th)|\D|^{1/2}\|_{\mathfrak{S}_{1,1}},\|\gamma_*+th\|_{\mathcal{Y}}\}\notag\\
    &+\frac{M(\alpha,c)}{c^2}\|T(\gamma_*+th)-\gamma_*-th\|_{X_c\bigcap \mathcal{Y}_c}\leq Kq+C_{\mathcal{Y}}+C_{M}q\frac{\alpha}{c^2}\leq R.
\end{align}
Thus, we know $\gamma_*+th\in \mathcal{U}_R$.

Now from Lemma \ref{lem:retra}, the limit $\theta(\gamma)$ exists for any $\gamma\in \mathcal{U}_R$. Hence $E(\gamma_*)$ and $E(\gamma_*+th)$ are well-defined. Then from \eqref{min-retra} and Theorem \ref{th:penal}, for any $t\in [-t_0,t_0]$,
\[
E(\gamma_*+th)\geq \min_{\gamma\in \mathcal{U}_R}E(\gamma)= E(\gamma_*).
\]
By definition of $\gamma_*$ and $h$, we also have
\[
\fint_{Q_{\ell}^*}\Tr_{L^2_\xi}[(\D_{\gamma_*,\xi}-\epsilon_P)h_{\xi}]\,d\xi=0.
\]
Then we deduce from Theorem \ref{approx} that
\[
 \frac{\alpha }{2}\widetilde{\Tr}_{L^2}[V_h h]+\alpha_c^2  Err(t,\gamma_*,h)\geq 0
\]
with $|Err(t,\gamma_*,h)|\leq (2+\frac{\epsilon_P}{c^2})\mathcal{N}_{\gamma_*+th}(h)$. This ends the proof.
\end{proof}
Let $B(\xi,\lambda)$ denote the ball of radius $\lambda >0 $ centered at $\xi \in \mathbb{R}^3$. As $|\omega|\neq 0$, we may find two points $\xi_1$ and $\xi_2$ in $\omega$ such that $|\omega\bigcap B(\xi_j,\lambda)|\neq 0$ for $j=1,2$. In particularly, for $\lambda$ small enough, we have $\omega\bigcap B(\xi_j,\lambda)=B(\xi_j,\lambda)$.  As in \cite{ghimenti2009properties}, we introduce the operator $h^{\lambda}$ defined by
\[
h^\lambda_\xi=\eta_\lambda(\xi)\left|\psi(\xi)\right>\left<\psi(\xi)\right|-\eta_\lambda'(\xi)\left|\psi'(\xi)\right>\left<\psi'(\xi)\right|
\]
for every $\xi\in  Q_\ell^*$, where
\[
\eta_\lambda=\frac{\mathbbm{1}_{\omega\bigcap B(\xi_1,\lambda)}}{|\omega\bigcap B(\xi_1,\lambda)|},\quad \eta_\lambda'=\frac{\mathbbm{1}_{\omega\bigcap B(\xi_2,\lambda)}}{|\omega\bigcap B(\xi_2,\lambda)|}.
\]
As shown in the proof of \cite[Lemma 9]{ghimenti2009properties}, we have
\begin{align}\label{eq:4.10}
    \widetilde{\Tr}_{L^2}[V_{h^\lambda} h^\lambda]\leq C-\frac{4\pi}{\lambda^2}.
\end{align}
On the other hand, 
\begin{lemma}\label{lem:contra2}
There exists a constant $C_{\rm cri}(q,z)>0$ large enough and only depending on $z$ and $q$ such that for any $\alpha_c\leq \frac{4\pi}{C_{\rm cri}(q,z)}$, $c\geq 2(C_Gz+C_{EE}q)$ and any $\lambda$ small enough, 
\[
\overline{\lim_{t\to 0}}|Err(t,\gamma_*,h^\lambda)|\leq (1+\frac{\epsilon_{P}}{c^2})\mathcal{N}_{\gamma_*}(h^\lambda)\leq C_{\rm cri}(q,z)+\frac{C_{\rm cri}(q,z)}{\lambda^2}.
\]
\end{lemma}
\begin{proof}
In this proof, $C$ denotes several different positive constant independent of $\alpha$ and $c$.

Let $R$ be given as in Lemma \ref{lem:contra0}. Then we choose first $C_{\rm cri}(q,z)\geq 16\pi C_{EE}R$. Then we have $\alpha_c\leq \frac{4\pi}{C_{\rm cri}(q,z)} \leq \frac{1}{4C_{EE}R}$. Thus Lemma \ref{lem:contra0} holds, and $\kappa(\alpha,c)\leq \frac{1}{2}$ and $L(\alpha,c)\leq \frac{1}{2}$. 
As a result, from Lemma \ref{lem:contra0}, we infer $\overline{\lim}_{t\to 0}|Err(t,\gamma_*,h^\lambda)|\leq (2+\frac{\epsilon_P}{c^2})\mathcal{N}_{\gamma_*}(h^\lambda)$ for any $\lambda>0$. In addition, according to Theorem \ref{th:penal} and Lemma \ref{lem:prop-spec}, there exists a constant $C$ independent of $\alpha$ and $c$ such that $(2+\frac{\epsilon_P}{c^2})\leq C$. Hence,
\begin{align*}
    \overline{\lim}_{t\to 0}|Err(t,\gamma_*,h^\lambda)|\leq (2+\frac{\epsilon_P}{c^2})\mathcal{N}_{\gamma_*}(h^\lambda)\leq C\mathcal{N}_{\gamma_*}(h^\lambda).
\end{align*}

We then prove that the terms appearing in $\mathcal{N}_{\gamma_*}(h^\lambda)$ satisfy respectively
\begin{align}
    \frac{1}{c}\fint_{Q_{\ell}^*}\|h^\lambda\|_{X\bigcap \mathcal{Y}(\xi)}\|\gamma_{*,\xi}|\D_\xi|^{1/2}\|_{\mathfrak{S}_1(\xi)}d\xi\leq C,\label{eq:lem4.7-1}\\
    \fint_{Q_{\ell}^*}\|h^\lambda\|_{X\bigcap \mathcal{Y}(\xi)}^2 \|\gamma_{*,\xi}\|_{\mathfrak{S}_1(\xi)}\,d\xi\leq C+\frac{C}{\lambda^2},\label{eq:lem4.7-2}
\end{align}
and
\begin{align}\label{eq:lem4.7-3}
    \sup_{\xi\in Q_\ell^*}\fint_{Q_{\ell}^*}\frac{\|h^\lambda\|_{X\bigcap \mathcal{Y}(\xi')}}{|\xi-\xi'|^2}\,d\xi'\leq C+\frac{C}{\lambda},
\end{align}
Once the above estimates are established, then there exists $C'_{\rm cri}(q,z)$ independent of $\alpha$ and $c$ such that $\mathcal{N}_{\gamma_*}(h^\lambda)\leq C'_{\rm cri}(q,z)+\frac{C'_{\rm cri}(q,z)}{\lambda^2}$ since $\frac{1}{2}\leq 1-L(\alpha,c)\leq 1$ and $\frac{1}{2}\leq 1-\kappa(\alpha,c)\leq \lambda_0(\alpha,c)\leq 1$ independent of $\alpha$ and $c$. Hence this lemma follows by choosing 
\begin{align*}
    C_{\rm cri}(q,z)=\max\{16\pi C_{EE}R,C'_{\rm cri}(q,z)\}.
\end{align*}

According to Theorem \ref{th:penal}, for any $\xi\in Q_\ell^*$, $0\leq \gamma_{*,\xi}\leq \mathbbm{1}_{(0,\nu]}(\D_{\gamma_*,\xi})$. As $\nu\leq c^*(q+1)$, by Lemma \ref{lem:prop-spec}, we have $\rank(\gamma_{*,\xi})\leq M+q$. Thus,
\begin{align}\label{eq:gamma-sigma1}
    \|\gamma_{*,\xi}\|_{\mathfrak{S}_1(\xi)}\leq \rank(\gamma_{*,\xi})\leq q+M.
\end{align}
On the other hand, we write $
\gamma_{*,\xi}=\sum_{k=1}^{M+q}\mu_k(\xi)\left|\psi_k(\xi)\right>\left<\psi_k(\xi)\right|$
where $0\leq \mu_k(\xi)\leq 1$ and $\psi_k(\xi)$ is the eigenfunction of $\D_{\gamma_*,\xi}$ with eigenvalue $0\leq \lambda_k(\xi)\leq c^*(q+1)$. By Lemma \ref{lem:unibound}, we finally get
\[
\frac{1}{c}\|\gamma_{*,\xi}|\D_\xi|^{1/2}\|_{\mathfrak{S}_1(\xi)}\leq \|\gamma_{*,\xi}(1-\Delta_\xi)^{1/4}\|_{\mathfrak{S}_1(\xi)}\leq \sum_{k=1}^{q+M}\|\psi_k(\xi)\|_{H^1_{\xi}(Q_\ell)}\leq C.
\]
It is easy to see that
\[
\|h^\lambda\|_X\leq |Q_\ell^*|^{-1}\fint_{\omega\bigcap B(\xi_1,\lambda)}\|\psi(\xi)\|_{H^{1/2}_{\rm per}}^2d\xi+|Q_\ell^*|^{-1}\fint_{\omega\bigcap B(\xi_2,\lambda)}\|\psi'(\xi)\|_{H^{1/2}_{\rm per}}^2d\xi.
\]
By Lemma \ref{lem:unibound} again, we get
\begin{align}\label{eq:h-X}
    \|h^\lambda\|_X \leq C.
\end{align}
We now turn to the study of the term $\|h^\lambda\|_{\mathcal{Y}(\xi)}$. First, 
\[
\|h^\lambda\|_{\mathcal{Y}(\xi)}=\fint_{Q_{\ell}^*}\frac{\|h_{\xi'}^\lambda\|_{\mathcal{B}(L^2_{\xi'})}}{|\xi-\xi'|^2}\,d\xi'\leq |Q_\ell^*|^{-1}\sum_{k=1,2}\fint_{\omega\bigcap B(\xi_k,\lambda)}\frac{d\xi'}{|\xi-\xi'|^2}.
\]
Thus,
\[
\frac{1}{c}\fint_{Q_{\ell}^*}\|h^\lambda\|_{X\bigcap \mathcal{Y}(\xi)}\|\gamma_{*,\xi}|\D_\xi|^{1/2}\|_{\mathfrak{S}_1(\xi)}d\xi\leq C+C\sum_{k=1,2}\fint_{\omega\bigcap B(\xi_k,\lambda)}\fint_{Q_\ell^*}\frac{1}{|\xi-\xi'|^2}d\xi d\xi'.
\]
As $\sup_{\xi'\in Q_\ell^*}\int_{Q_\ell^*}\frac{1}{|\xi-\xi'|^2}d\xi<+\infty$ independent of $\lambda$, we finally get
\[
\frac{1}{c}\fint_{Q_{\ell}^*}\|h^\lambda\|_{X\bigcap \mathcal{Y}(\xi)}\|\gamma_{*,\xi}|\D_\xi|^{1/2}\|_{\mathfrak{S}_1(\xi)}d\xi\leq C.
\]
This gives \eqref{eq:lem4.7-1}. Analogously, using \eqref{eq:gamma-sigma1},
\begin{align*}
    \fint_{Q_\ell^*}\|h^\lambda\|_{X\bigcap \mathcal{Y}(\xi)}\|\gamma_{*,\xi}\|_{\mathfrak{S}_1(\xi)}d\xi\leq C.
\end{align*}
Then, for Eq. \eqref{eq:lem4.7-2}, for $\lambda$ small enough, we have
\begin{align*}
    \fint_{Q_{\ell}^*}\|h^\lambda\|_{X\bigcap \mathcal{Y}(\xi)}^2 \|\gamma_{*,\xi}\|_{\mathfrak{S}_1(\xi)}d\xi&\leq C\|h^\lambda\|_{X\bigcap \mathcal{Y}}\fint_{Q_\ell^*}\|h^\lambda\|_{X\bigcap \mathcal{Y}(\xi)}\|\gamma_{*,\xi}\|_{\mathfrak{S}_1(\xi)}d\xi\\
&\leq C+ C\sum_{k=1,2}\sup_{\xi\in Q_\ell^*}\fint_{B(\xi_k,\lambda)}\frac{1}{|\xi-\xi'|^2}d\xi'\leq C+\frac{C}{\lambda^2}.
\end{align*}
We turn now to the last estimate \eqref{eq:lem4.7-3}. From \eqref{eq:h-X}, we have
\[
\sup_{\xi\in Q_\ell^*}\fint_{Q_{\ell}^*}\frac{\|h^\lambda\|_{X}}{|\xi-\xi'|^2}\,d\xi'\leq C.
\]
On the other hand,  according to \cite[Chapter 5.10, Formula (3)]{lieb2001analysis}, 
\begin{align*}
    \int_{\mathbb{R}^3}\frac{d\xi'}{|\xi-\xi'|^2|\xi'-\xi''|^2} = \frac{\pi^3}{|\xi-\xi''|}.
\end{align*}
Hence, for $\lambda$ small enough,
\begin{align*}
 \sup_{\xi\in Q_\ell^*}\fint_{Q_{\ell}^*}\frac{\|h^\lambda\|_{\mathcal{Y}(\xi')}}{|\xi-\xi'|^2}\,d\xi'&\leq C\sum_{k=1,2}\sup_{\xi\in Q_\ell^*}\fint_{\omega\bigcap B(\xi_k,\lambda)}d\xi''\int_{Q_\ell^*}\frac{d\xi'}{|\xi-\xi'|^2|\xi'-\xi''|^2}\\
    &\leq C\sum_{k=1,2}\sup_{\xi\in Q_\ell^*}\fint_{\omega\bigcap B(\xi_k,\lambda)}\frac{d\xi''}{|\xi-\xi''|}\leq C\sum_{k=1,2}\fint_{ B(\xi_k,\lambda)}\frac{d\xi''}{|\xi''-\xi_k|}\leq \frac{C}{\lambda}.
\end{align*}
Thus we deduce 
\[
\sup_{\xi\in Q_\ell^*}\fint_{Q_{\ell}^*}\frac{\|h^\lambda\|_{X\bigcap \mathcal{Y}(\xi')}}{|\xi-\xi'|^2}d\xi'\leq \frac{C}{\lambda}.
\]
This ends the proof of Lemma~\ref{lem:contra2}.
\end{proof}
This ends the proof of Theorem \ref{th:prop}. Indeed, by using \eqref{eq:4.10} and Lemma \ref{lem:contra2}, under Assumption \ref{ass:2} we get
\begin{align*}
    \overline{\lim}_{t\to 0}\left(\frac{\alpha }{2}\widetilde{\Tr}_{L^2}[V_{h^\lambda} h^\lambda]+\alpha_c^2  Err(t,\gamma_*,h^\lambda)\right)\leq C+C_{\rm cri}(q,z) +\frac{C_{\rm cri}(q,z)\alpha_c^2-4\pi\alpha}{\lambda^2}. 
\end{align*}
 Thus, for $\alpha>0$ and $\alpha_c$ small enough satisfying \begin{align}\label{eq:C_cri}
\alpha_c<\frac{4\pi}{C_{\rm cri}(q,z)},
\end{align}
we have 
\begin{align*}
   \lim_{\lambda\to 0} \left(\frac{\alpha }{2}\widetilde{\Tr}_{L^2}[V_{h^\lambda} h^\lambda]+\alpha_c^2  Error_\gamma(h^\lambda)\right)\leq C+C_{\rm cri}(q,z)+\lim_{\lambda\to 0}\frac{(C_{\rm cri}(q,z)\alpha_c-4\pi)\alpha}{\lambda^2} =-\infty.
\end{align*}
We reach a contradiction with \eqref{eq:4.9} whenever Assumption \ref{ass:2} is satisfied.  Hence the proof of Theorem \ref{th:prop} is completed.
\section{Existence of the retraction}\label{sec:retra}
In this section, we argue as in \cite[Proposition 5.8]{crystals} and \cite{Ser09} to give a sketch of the proof of Lemma \ref{lem:retra}. It is based on the following result.
\begin{lemma} 
Let $\kappa(\alpha,c)<1$. Recall that $A(\alpha,c):=\frac{\alpha_c}{2}C_{EE}(1-\kappa(\alpha,c))^{-1/2}\lambda_0(\alpha,c)^{-1/2}$. Then for any $\gamma,\gamma'\in \Gamma_q$ and $\xi \in Q_\ell^*$
    \begin{align}\label{eq:P-P}
        \||\D_\xi|^{1/2}(P^+_{\gamma,\xi}-P^+_{\gamma',\xi})\|_{\mathcal{B}(L^2_\xi)}\leq \frac{A(\alpha,c)}{(1+C_\mathcal{Y})} \|\gamma-\gamma'\|_{X\bigcap \mathcal{Y}(\xi)}\leq \frac{A(\alpha,c)}{c(1+C_{\mathcal{Y}})}\|\gamma-\gamma'\|_{X_c\bigcap \mathcal{Y}_c(\xi)}
    \end{align}
and
\begin{align}\label{eq:T2-T}
    \|T^2(\gamma)-T(\gamma)\|_{X_c\bigcap \mathcal{Y}_c}&\leq 2A(\alpha,c)\big(\frac{1}{c} \max\{\|T(\gamma)|\D|^{1/2}\|_{\mathfrak{S}_{1,1}},\|\gamma\|_{\mathcal{Y}}\}\big.\notag\\
    &\qquad\big.+\frac{A(\alpha,c)q}{c^2}\|T(\gamma)-\gamma\|_{X_c\bigcap \mathcal{Y}_c}\big)\|\gamma-\gamma'\|_{X_c\bigcap \mathcal{Y}_c}.
\end{align}
\end{lemma}
\begin{proof}
By Taylor's formula, we have
\begin{align}\label{eq:Taylor}
    P^\pm_{\gamma}=\frac{1}{2}\pm \frac{1}{2\pi}\int_{-\infty}^{+\infty}(\D_{\gamma}-iz)^{-1}dz.
\end{align}
Thus,
\[
P^\pm_{\gamma}-P^\pm_{\gamma'}=\pm\frac{\alpha}{2\pi}\int_{-\infty}^{+\infty}(\D_{\gamma}-iz)^{-1}V_{\gamma'-\gamma}(\D_{\gamma'}-iz)^{-1}dz.
\]
Now according to \eqref{3.eq:5.1'}, \eqref{eq:D-D} and \eqref{eq:A.9} and the formula
\begin{align}\label{eq:arctan}
    \int_{-\infty}^{+\infty}\frac{|B|}{|B|^2+z^2}dz=\pi,\quad \textrm{with } B\neq 0,
\end{align}
for any $\psi_\xi,\phi_\xi\in L^2_\xi$, we obtain
\begin{align*}
    \MoveEqLeft \left(\psi_{\xi},|\D_{\xi}|^{1/2}(P^+_{\gamma,\xi}-P^+_{\gamma',\xi})\phi_{\xi}\right)_{L^2_{\xi}}\\
    &= \frac{\alpha}{2\pi}\int_{-\infty}^{+\infty}\left(\psi_{\xi},|\D_{\xi}|^{1/2}(\D_{\gamma,\xi}-iz)^{-1}V_{\gamma-\gamma',\xi}(\D_{\gamma',\xi}-iz)^{-1} \phi_{\xi}\right)_{L^2_{\xi}}dz\\
    &\leq \frac{\alpha}{2\pi}\|V_{\gamma-\gamma',\xi}\|_{\mathcal{B}(L^2_{\xi})}\left(\int_{-\infty}^{+\infty}\|(\D_{\gamma,\xi}-iz)^{-1}|\D_{\xi}|^{1/2}\psi_{\xi}\|_{L^2_{\xi}}^2dz\right)^{1/2}\\
    &\quad\times\left(\int_{-\infty}^{+\infty}\|(\D_{\gamma',\xi}-iz)^{-1}\phi_{\xi}\|_{L^2_{\xi}}^2dz\right)^{1/2}\\
    &\leq \frac{\alpha}{2}\|V_{\gamma-\gamma'}\|_{Y}\||\D|^{1/2}|\D_{\gamma}|^{-1/2}\|_Y\||\D_{\gamma'}|^{-1/2}\|_Y\|\psi_{\xi}\|_{L^2_{\xi}}\|\phi\|_{L^2_{\xi}}\\
    &\leq \frac{C_{EE}}{2(1+C_\mathcal{Y})}(1-\kappa(\alpha,c))^{-1/2}\lambda_0(\alpha,c)^{-1/2}\alpha_c \|\gamma-\gamma'\|_{X\bigcap \mathcal{Y}}\|\psi\|_{\mathcal{H}}\|\phi\|_{\mathcal{H}}\\
    &\leq \frac{C_{EE}}{2(1+C_{\mathcal{Y}})}(1-\kappa(\alpha,c))^{-1/2}\lambda_0(\alpha,c)^{-1/2}\frac{\alpha_c}{c}\|\gamma-\gamma'\|_{X_c\bigcap \mathcal{Y}_c}\|\psi\|_{\mathcal{H}}\|\phi\|_{\mathcal{H}}.
\end{align*}
Hence \eqref{eq:P-P}.  Now, we turn to the proof of \eqref{eq:T2-T}. We have
\begin{align*}
    \MoveEqLeft T^2(\gamma)-T(\gamma)=(P^+_{T(\gamma)}-P^+_{\gamma})T(\gamma)P^+_{T(\gamma)}\\
    &+T(\gamma)(P^+_{T(\gamma)}-P^+_{\gamma})+(P^+_{T(\gamma)}-P^+_{\gamma})T(\gamma)(P^+_{T(\gamma)}-P^+_{\gamma}).
\end{align*}
Hence,
\begin{align*}
    \MoveEqLeft \|T^2(\gamma)-T(\gamma)\|_{X_c\bigcap \mathcal{Y}_c}\leq 2\|(P^+_{T(\gamma)}-P^+_{\gamma})T(\gamma)\|_{X_c\bigcap \mathcal{Y}_c}+\|(P^+_{T(\gamma)}-P^+_{\gamma})T(\gamma)(P^+_{T(\gamma)}-P^+_{\gamma})\|_{X_c\bigcap \mathcal{Y}_c}.
\end{align*}
Since $c^2\leq |\D|$, we have
\begin{align*}
\|T(\gamma)(P^+_{T(\gamma)}-P^+_\gamma)\|_{X_c\bigcap \mathcal{Y}_c}\leq \||\D|^{1/2}(P^+_{T(\gamma)}-P^+_{\gamma})\|_Y\max\{\|T(\gamma)|\D|^{1/2}\|_{\mathfrak{S}_{1,1}},\|T(\gamma)\|_{\mathcal{Y}}\},
\end{align*}
and
\begin{align*}
\MoveEqLeft\|(P^+_{T(\gamma)}-P^+_\gamma)T(\gamma)(P^+_{T(\gamma)}-P^+_\gamma)\|_{X_c\bigcap \mathcal{Y}_c}\leq \||\D|^{1/2}(P^+_{T(\gamma)}-P^+_{\gamma})\|_Y^2\|T(\gamma)\|_{\mathfrak{S}_{1,1}\bigcap \mathcal{Y}}.
\end{align*}
Observe from \eqref{eq:T-Y} that $\|T(\gamma)\|_{\mathfrak{S}_{1,1}}\leq \|\gamma\|_{\mathfrak{S}_{1,1}}\leq q$ and $\|T(\gamma)\|_{\mathcal{Y}}\leq \|\gamma\|_{\mathcal{Y}}\leq C_\mathcal{Y}$. Then we have $\|T(\gamma)\|_{\mathfrak{S}_{1,1}\bigcap \mathcal{Y}} \leq (1+C_\mathcal{Y})q$. Then using \eqref{eq:P-P}, we obtain
\begin{align*}
    \|T^2(\gamma)-T(\gamma)\|_{X_c\bigcap Y_c}\leq& 2A(\alpha,c)\left(\frac{1}{c}\max\{\|T(\gamma)|\D|^{1/2}\|_{\mathfrak{S}_{1,1}},\|T(\gamma)\|_{\mathcal{Y}}\}\right.\\
    &\left.+\frac{A(\alpha,c)q}{2c^2}\|\gamma-T(\gamma)\|_{X_c\bigcap \mathcal{Y}_c}\right)\|\gamma-T(\gamma)\|_{X_c\bigcap \mathcal{Y}_c}.
\end{align*}
\end{proof}
We turn now to the following.
\begin{proof}[Proof of Lemma \ref{lem:retra}]
First of all, observe that $T(\gamma)\leq \gamma\leq \mathbbm{1}_{L^2(\mathbb{R}^3)}$, $\|T(\gamma)\|_{\mathfrak{S}_{1,1}}\leq \|\gamma\|_{\mathfrak{S}_{1,1}}\leq q$, $\|T(\gamma)\|_{\mathcal{Y}}\leq \|\gamma\|_{\mathcal{Y}}$, and from \eqref{eq:DPD}
\[
\|T(\gamma)\|_{X}\leq \frac{1-\kappa(\alpha,c)}{1+\kappa(\alpha,c)}\|\gamma\|_X.
\]
Then $T(\gamma)\in \Gamma_{\leq q}$. We are going to prove that $T$ maps $\mathcal{U}_R$ into $\mathcal{U}_R$. For $\gamma\in\mathcal{U}_R$, we have
\begin{align}\label{eq:5.4}
    \MoveEqLeft \|T(\gamma)|\D|^{1/2}\|_{\mathfrak{S}_{1,1}}\leq \|\gamma|\D|^{1/2}\|_{\mathfrak{S}_{1,1}}+\frac{1}{c}\|\gamma-T(\gamma)\|_{X_c}
\end{align}
As $M(\alpha,c)\geq \frac{1+A(\alpha,c)q}{2}$, \eqref{eq:T2-T} implies 
\begin{align}\label{eq:5.5}
    \|T^2(\gamma)-T(\gamma)\|_{X_c}\leq L(\alpha,c)\|T(\gamma)-\gamma\|_{X_c},
\end{align}
with $L(\alpha,c)=2A(\alpha,c)R$, and
\begin{align}\label{eq:5.6}
    \frac{1}{c}\|T(\gamma)\|_{\mathcal{Y}}\leq \frac{1}{c}\|\gamma\|_{\mathcal{Y}}+\frac{1}{c^2}\|T(\gamma)-\gamma\|_{\mathcal{Y}_c}.
\end{align}
Moreover, as $M(\alpha,c)\geq \frac{1}{1-2A(\alpha,c)R}$, we have $1+M(\alpha,c)L(\alpha,c)\leq M$. Then from \eqref{eq:5.4}-\eqref{eq:5.6} and the fact that $\gamma\in \mathcal{U}_R$,
\begin{align*}
    \MoveEqLeft \frac{1}{c}\max\{\|T(\gamma)|\D|^{1/2}\|_{\mathfrak{S}_{1,1}},\|T(\gamma)\|_{\mathcal{Y}}\}+\frac{M(\alpha,c)}{c^2}\|T^2(\gamma)-T(\gamma)\|_{X_c\bigcap \mathcal{Y}_c}\\
    &\qquad\leq \frac{1}{c}\max\{\|\gamma|\D|^{1/2}\|_{\mathfrak{S}_{1,1}},\|\gamma\|_{\mathcal{Y}}\}+\frac{1+M(\alpha,c)L(\alpha,c)}{c^2}\|\gamma-T(\gamma)\|_{X_c\bigcap \mathcal{Y}_c}<R.
\end{align*}
Then $T(\gamma)\in \mathcal{U}_R$. Thus $T$ maps $\mathcal{U}_R$ into $\mathcal{U}_R$, so does $T^n$ for any $n\in\mathbb{N}$. 

Next, we prove the existence of $\theta$. From \eqref{eq:T2-T}, we know 
\begin{align}\label{eq:retra'}
    \|T^n(\gamma)-T^{n-1}(\gamma)\|_{X_c\bigcap \mathcal{Y}_c}\leq L(\alpha,c)\|T^{n-1}(\gamma)-T^{n-2}(\gamma)\|_{X_c\bigcap \mathcal{Y}_c}
\end{align}
and
\[
\|T^n(\gamma)\|_{X_c\bigcap \mathcal{Y}_c}\leq \sum_{i=1}^n\|T^i(\gamma)-T^{i-1}(\gamma)\|_{X_c\bigcap \mathcal{Y}_c}+\|\gamma\|_{X}\leq \frac{1}{1-L(\alpha,c)}\|T(\gamma)-\gamma\|_{X_c\bigcap \mathcal{Y}_c}+\|\gamma\|_{X_c\bigcap \mathcal{Y}_c}.
\]
This implies that $\|T^n(\gamma)\|_{X_c\bigcap \mathcal{Y}_c}$ is uniformly bounded with respect to $n\in\mathbb{N}$, and $(T^n(\gamma))_n$ is a Cauchy sequence in $X_c\bigcap \mathcal{Y}_c$. Thus for any $\gamma\in\mathcal{U}_R$, the retraction $\theta(\gamma):=\lim_{n\to +\infty}T^n(\gamma)$ exists in $X_c\bigcap \mathcal{Y}_c$. Furthermore, we have
\[
\|\theta(\gamma)-T^n(\gamma)\|_{X_c\bigcap \mathcal{Y}_c}\leq \frac{L(\alpha,c)^n}{1-L(\alpha,c)}\|T(\gamma)-\gamma\|_{X_c\bigcap \mathcal{Y}_c}.
\]
From \eqref{eq:P-P}, it can be deduced directly that $T$ is continuous on $\mathcal{U}_R$. Finally, the fact that $\theta(\gamma)\in\Gamma_q^+$ for any $\gamma\in \mathcal{U}_R$ follows from the fact 
\begin{align*}
 \|T(\theta(\gamma))-\theta(\gamma)\|_{X_c\bigcap \mathcal{Y}_c}&= \lim_{n\to +\infty}\|T^{n+1}(\gamma)-T^{n}(\gamma)\|_{X_c\bigcap \mathcal{Y}_c}= 0.
\end{align*}
This ends the proof.
\end{proof}

\section{Second-order expansion of \texorpdfstring{$E(\gamma)$}{}}\label{sec:5}
This section is devoted to the proof of Theorem \ref{approx}. We fix $\alpha$ and $c$. For future convenience, we denote $L:=L(\alpha,c)$, $\kappa:=\kappa(\alpha,c)$ and $\lambda_0:=\lambda_0(\alpha,c)$. The main ingredient is the following which is essentially the same as in \cite{meng}.
\begin{proposition}\label{main theorem}
Let $R>1$ be fixed, and let $\kappa<1$ and $L<1$ be given as in Lemma \ref{lem:retra}. For any $\gamma\in\mathcal{U}_R$ and any $g\in\Gamma_{\leq q}$, if $P^+_{g}\gamma P^+_{g}=\gamma$, we have
\begin{align}\label{E-E}
    |E(\gamma)-(\mathcal{E}(\gamma)-\epsilon_P\widetilde{\Tr}_{L^2}(\gamma))|\leq (2+\frac{\epsilon_P}{c^2})\alpha_c^2 \mathcal{N}_\gamma(\gamma-g),
\end{align}
where $\mathcal{N}_\gamma(h)$ is given by \eqref{eq:error}, i.e.,
\begin{align*}
\mathcal{N}_\gamma(h)&=\frac{ C_{EE}^2}{2(1-\kappa)^{2}\lambda_0}\fint_{Q_{\ell}^*}\|h\|_{X\bigcap \mathcal{Y}(\xi)}^2 \|\gamma_\xi\|_{\mathfrak{S}_1(\xi)}d\xi\\
 &+(q\alpha_c^2+R\alpha_c^2 +\alpha_c) \frac{10C_{EE}^4}{(1-\kappa)^4\lambda_0^{5/2}(1-L)^2} \\
  &\qquad \times\left(\frac{1}{c}\fint_{Q_{\ell}^*}\|h\|_{X\bigcap \mathcal{Y}(\xi)}\|\gamma_\xi|\D_\xi|^{1/2}\|_{\mathfrak{S}_1(\xi)}d\xi
+\frac{1}{c}\sup_{\xi\in Q_\ell^*}\fint_{Q_{\ell}^*}\frac{\|h\|_{X\bigcap \mathcal{Y}(\xi')}}{|\xi-\xi'|^2}d\xi'\right)^2.
\end{align*}
\end{proposition}
\begin{remark}
Actually, for atoms and molecules, using a finer version of this estimate, it is shown in \cite{meng} that the DF model is an approximation of the electron-positron Hartree-Fock model (see, e.g., \cite{barbaroux2005hartree} for this model).
\end{remark}
We first use this proposition to prove Theorem \ref{approx}.
\begin{proof}[Proof of Theorem \ref{approx}]
Let 
\[
Err(t,\gamma,h)=\frac{1}{\alpha_c^2 t^2}(E(\gamma+th)-\mathcal{E}(\gamma+th)).
\]
Since $\gamma\in \Gamma_{\leq q}^+$, we have
\[
E(\gamma+th)=E(\gamma)+\mathcal{E}(\gamma+th)-\mathcal{E}(\gamma)+t^2\alpha_c^2 Err(t,\gamma,h).
\]
This is \eqref{eq:second}. Then replacing $g,\gamma$ by $\gamma,\gamma+th$ respectively in Proposition \ref{main theorem}, we finally get
\begin{align}
    |Err(t,\gamma,h)|\leq (2+\frac{\epsilon_P}{c^2})t^{-2}\mathcal{N}_{\gamma+th}(th)=(2+\frac{\epsilon_P}{c^2})\mathcal{N}_{\gamma+th}(h).
\end{align}
\end{proof}

\subsection{Proof of Proposition \ref{main theorem}}\label{sec:5.2}
We first consider the error bound between $E(\gamma)$ and $\mathcal{E}(\gamma)$ for any $\gamma\in \mathcal{U}_R$.
\begin{lemma}\label{lem:6.3}
    Let $R\geq 1$ be fixed. Assume that $\kappa<1$ and $L<1$ as in Lemma \ref{lem:retra}. Let $C_{\kappa,L}:=\frac{5C_{EE}^2}{(1-\kappa)^2\lambda_0^{3/2}(1-L)^2}$. Then for any $\gamma\in \mathcal{U}_R$,
    \begin{align}\label{eq:6.3}
        |E(\gamma)-\mathcal{E}(\gamma)|&\leq \frac{2c^2+\epsilon_P}{c^2}\left[C_{\kappa,L}(R\alpha_c+q\alpha_c +1)\frac{\alpha_c}{c^2}\|T(\gamma)-\gamma\|_{X_c\bigcap \mathcal{Y}_c}^2+\|P^-_{\gamma}\gamma P^-_{\gamma}\|_{X_c}\right].
    \end{align}
\end{lemma}
This is an immediate result of the following.
\begin{lemma}\label{lem:t-T}
Let $\kappa\leq 1$ and $L<1$. For any $\gamma\in \mathcal{U}_R$,
\begin{align}\label{theta-T+}
    \|P^{+}_{\gamma}(\theta(\gamma)-T(\gamma))P^{+}_{\gamma}\|_{X_c\bigcap \mathcal{Y}_c} \leq C_{\kappa,L} R\frac{\alpha_c^2}{c^2}\|T(\gamma)-\gamma\|_{X_c\bigcap \mathcal{Y}_c}^2
\end{align}
and
\begin{align}\label{theta-T-} \|P^-_{\gamma}\theta(\gamma)P^-_{\gamma}\|_{X_c\bigcap \mathcal{Y}_c}\leq C_{\kappa,L}q\frac{\alpha_c^2}{c^2}\|T(\gamma)-\gamma\|_{X_c\bigcap \mathcal{Y}_c}^2.
\end{align}
\end{lemma}
We first use it to prove Lemma \ref{lem:6.3} and we postpone the proof of \ref{lem:t-T} until Section \ref{sec:6.3}.
\begin{proof}[Proof of Lemma \ref{lem:6.3}]
    Notice that
    \begin{align}\label{eq:E-E}
        E(\gamma)-\mathcal{E}(\gamma)=\widetilde{\Tr}_{L^2} [\D_{\gamma}(\theta(\gamma)-\gamma)]+\frac{\alpha}{2} \widetilde{\Tr}_{L^2}[V_{\theta(\gamma)-\gamma}(\theta(\gamma)-\gamma)]-\epsilon_P\widetilde{\Tr}_{L^2}[\theta(\gamma)- \gamma].
    \end{align}
To end the proof, we calculate each term on the right-hand side separately. 

\textbf{Estimate on $\widetilde{\Tr}_{L^2}[\D_{\gamma}(\theta(\gamma)-\gamma)]$.} We consider the first term on the right-hand side of \eqref{eq:E-E}. Notice that $T(\gamma)=P^+_{\gamma}\gamma P^+_{\gamma}$. We have
\begin{align*}
   \widetilde{\Tr}_{L^2}[\D_{\gamma}(\theta(\gamma)-\gamma)]&=\widetilde{\Tr}_{L^2}[\D_{\gamma} (P^+_{\gamma}+P^-_{\gamma})(\theta(\gamma)-\gamma)(P^+_{\gamma}+P^-_{\gamma})]\\
    &=\widetilde{\Tr}_{L^2}[|\D_{\gamma}|P^+_{\gamma}(\theta(\gamma)-T(\gamma))P^+_{\gamma}]-\widetilde{\Tr}_{L^2} [|\D_{\gamma}|P^-_{\gamma}(\theta(\gamma)-\gamma)P^-_{\gamma}].
\end{align*}
Thanks to \eqref{eq:D-D}, \eqref{theta-T+} and $\kappa<1$, we have
\begin{align*}
    \left|\widetilde{\Tr}_{L^2}[|\D_{\gamma}|P^+_{\gamma}(\theta(\gamma)-T(\gamma))P^+_{\gamma}]\right|&\leq \||\D_{\gamma}|^{1/2}P^+_{\gamma}(\theta(\gamma)-T(\gamma))P^+_{\gamma}|\D_{\gamma}|^{1/2}\|_{\mathfrak{S}_{1,1}}\\
    &\leq 2\|P^+_{\gamma}(\theta(\gamma)-T(\gamma))P^+_{\gamma}\|_{X_c}\leq 2C_{\kappa,L}R\frac{\alpha_c^2}{c^2}\|T(\gamma)-\gamma\|_{X_c\bigcap \mathcal{Y}_c}^2.
\end{align*}
On the other hand, from \eqref{eq:D-D} and \eqref{theta-T-}, we infer
\begin{align*}
     \left|\widetilde{\Tr}_{L^2} [|\D_{\gamma}|P^-_{\gamma}(\theta(\gamma)-\gamma)P^-_{\gamma}]\right|&\leq \||\D_{\gamma}|^{1/2}P^-_{\gamma}\theta(\gamma)P^-_{\gamma}|\D_{\gamma}|^{1/2}\|_{\mathfrak{S}_{1,1}}+\||\D_{\gamma}|^{1/2}P^-_{\gamma} \gamma P^-_{\gamma}|\D_{\gamma}|^{1/2}\|_{\mathfrak{S}_{1,1}}\\
     &\leq 2(\|P^-_{\gamma}\theta(\gamma)P^-_{\gamma}\|_{X_c}+\|P^-_{\gamma} \gamma P^-_{\gamma}\|_{X_c})\\
     &\leq 2C_{\kappa,L}q\frac{\alpha_c^2}{c^2}\|T(\gamma)-\gamma\|^2_{X_c\bigcap \mathcal{Y}_c}+2\|P^-_{\gamma} \gamma P^-_{\gamma}\|_{X_c}.
\end{align*}
Then we conclude that
\begin{align*}
     \left|\widetilde{\Tr}_{L^2} [\D_{\gamma}(\theta(\gamma)-\gamma)]\right|\leq 2C_{\kappa,L}(R+q)\frac{\alpha_c^2}{c^2}\|T(\gamma)-\gamma\|^2_{X_c\bigcap \mathcal{Y}_c}+2\|P^-_{\gamma} \gamma P^-_{\gamma}\|_{X_c\bigcap \mathcal{Y}_c}.
\end{align*}

\textbf{Estimate on $\epsilon_P\widetilde{\Tr}_{L^2}[\theta(\gamma)- \gamma]$.} This term can be treated analogously. Actually, as $c^2\leq |\D|$, we have
\begin{align*}
    \epsilon_P|\widetilde{\Tr}_{L^2}[\theta(\gamma)- \gamma]|&\leq \epsilon_P\left|\widetilde{\Tr}_{L^2}[P^+_{\gamma}(\theta(\gamma)- \gamma)P^+_{\gamma}] +\widetilde{\Tr}_{L^2}[P^-_{\gamma}(\theta(\gamma)- \gamma)P^-_{\gamma}] \right|\\
    &\leq \frac{\epsilon_P}{c^2}\left(\|P^+_{\gamma}(\theta(\gamma)-T(\gamma))P^+_{\gamma}\|_{X_c}+ \|P^-_{\gamma}\theta(\gamma)P^-_{\gamma}\|_{X_c}+\|P^-_{\gamma} \gamma P^-_{\gamma}\|_{X_c}\right).
\end{align*}
Then proceeding as for the term $\Tr [\D_{\gamma}(\theta(\gamma)-\gamma)]$, we obtain
\begin{align*}
    \epsilon_P|\widetilde{\Tr}_{L^2}[\theta(\gamma)- \gamma]|\leq \frac{\epsilon_P}{c^2}C_{\kappa,L}q\frac{\alpha_c^2}{c^2}\|T(\gamma)-\gamma\|^2_{X_c\bigcap \mathcal{Y}_c}+\frac{\epsilon_P}{c^2}\|P^-_{\gamma} \gamma P^-_{\gamma}\|_{X_c\bigcap \mathcal{Y}_c}.
\end{align*}

\textbf{Estimate on $\frac{\alpha}{2} \widetilde{\Tr}_{L^2}[V_{\theta(\gamma)-\gamma}(\theta(\gamma)-\gamma)]$.} Using \eqref{3.eq:5.1'} and \eqref{eq:3.4}, we infer
\begin{align*}
  \MoveEqLeft  \alpha|\widetilde{\Tr}_{L^2}[V_{\theta(\gamma)-\gamma}(\theta(\gamma)-\gamma)]|\leq \frac{C_{EE}}{(1+C_\mathcal{Y})} \frac{\alpha_c }{c}\|\theta(\gamma)-\gamma\|_{X_c\bigcap \mathcal{Y}_c}\|\theta(\gamma)-\gamma\|_{\mathfrak{S}_{1,1}} \\
    &\qquad\leq \frac{C_{EE}}{(1+C_\mathcal{Y})(1-L)^2}\frac{\alpha_c}{c^2}\|T(\gamma)-\gamma\|_{X_c\bigcap \mathcal{Y}_c}^2 \leq C_{\kappa,L}\frac{\alpha_c}{c}\|T(\gamma)-\gamma\|_{X_c\bigcap \mathcal{Y}_c}.
\end{align*}

\textbf{Conclusion.} We then deduce the following estimate.
\begin{align*}
    |E(\gamma)-\mathcal{E}(\gamma)|&\leq \frac{2c^2+\epsilon_P}{c^2}\left[C_{\kappa,L}(R\alpha_c+q\alpha_c+1)\frac{\alpha_c}{c^2}\|T(\gamma)-\gamma\|_{X_c\bigcap \mathcal{Y}_c}^2+\|P^-_{\gamma}\gamma P^-_{\gamma}\|_{X_c\bigcap \mathcal{Y}_c}\right].
\end{align*}

\end{proof}

Now we consider the term $T(\gamma)-\gamma$ and $P^-_{\gamma}\gamma P^-_{\gamma}$ under the condition $\gamma=P^+_{g}\gamma P^-_{g}$.
\begin{lemma}\label{lem:T-I}
Let $\kappa<1$. Provided $g\in\Gamma_q$ and $\gamma\in \mathcal{U}_R$, if $P^+_{g}\gamma P^+_{g}=\gamma$, we have
\[
\|T(\gamma)-\gamma\|_{X_c}\leq  \frac{\sqrt{2} C_{EE}}{\lambda_0^{1/2}(1-\kappa)}  \alpha_c\fint_{Q_{\ell}^*}\|\gamma-g\|_{X\bigcap \mathcal{Y}(\xi)}\|\gamma_\xi|\D_\xi|^{1/2}\|_{\mathfrak{S}_1(\xi)}d\xi,
\]
\[
\|T(\gamma)-\gamma\|_{\mathcal{Y}_c}\leq  \frac{\sqrt{2} C_{EE}}{\lambda_0^{1/2}(1-\kappa)}  \alpha_c \sup_{\xi}\fint_{Q_{\ell}^*}\frac{\|\gamma-g\|_{X\bigcap \mathcal{Y}(\xi')}}{|\xi-\xi'|^2}d\xi',
\]
and
\begin{align*}
    \|P^-_{\gamma}\gamma P^-_{\gamma}\|_{X_c}\leq \frac{ C_{EE}^2}{2(1-\kappa)^{2}\lambda_0}\alpha_c^2\fint_{Q_{\ell}^*}\|g-\gamma\|_{X\bigcap \mathcal{Y}(\xi)}^2 \|\gamma_\xi\|_{\mathfrak{S}_1(\xi)}d\xi.
\end{align*}
\end{lemma}
\begin{proof}
Indeed, we have
\[
T(\gamma)-\gamma=(P^+_{\gamma}-P^+_{g})\gamma P^+_{\gamma}+P^+_g \gamma (P^+_{\gamma}-P^+_g).
\]
Then according to \eqref{eq:P-P}, \eqref{eq:DPD} and as $0<\kappa< 1$,
\begin{align*}
\|T_\xi(\gamma)-\gamma_\xi\|_{X_c(\xi)}&\leq \frac{2(1+\kappa)^{1/2}}{(1-\kappa)^{1/2}}\||\D_\xi|^{1/2}(P^+_{\gamma,\xi}-P^+_{g,\xi})\|_{\mathcal{B}(L^2_\xi)}\|\gamma_\xi|\D_\xi|^{1/2}\|_{\mathfrak{S}_1(\xi)}\\
&\leq \frac{\sqrt{2} C_{EE}}{(1+C_\mathcal{Y})\lambda_0^{1/2}(1-\kappa)}  \alpha_c\|\gamma-g\|_{X\bigcap \mathcal{Y}(\xi)}\|\gamma_\xi|\D_\xi|^{1/2}\|_{\mathfrak{S}_1(\xi)}.
\end{align*}
Thus,
\[
\|T(\gamma)-\gamma\|_{X_c}\leq \frac{\sqrt{2} C_{EE}}{\lambda_0^{1/2}(1-\kappa)}  \alpha_c\fint_{Q_{\ell}^*}\|\gamma-g\|_{X\bigcap \mathcal{Y}(\xi)}\|\gamma_\xi|\D_\xi|^{1/2}\|_{\mathfrak{S}_1(\xi)}d\xi.
\]
Analogously for the second estimate, the bound  $\|\gamma_\xi\|_{\mathcal{B}(L^2_\xi)}\leq 1$, $c^2\leq |\D|$, and \eqref{eq:T-Y} yield
\begin{align*}
   \|T(\gamma)-\gamma\|_{\mathcal{Y}_c} &\leq \sup_{\xi\in Q_\ell^*}\fint_{Q_{\ell}^*}\frac{2}{|\xi-\xi'|^2}\||\D_{\xi'}|^{1/2}(P^+_{\gamma,\xi'}-P^+_{g,\xi'})\|_{\mathcal{B}(L^2_{\xi'})}\|\gamma_{\xi'}\|_{\mathcal{B}(L^2_{\xi'})}d\xi'\\
&\leq \frac{\sqrt{2} C_{EE}}{\lambda_0^{1/2}(1-\kappa)}  \alpha_c \sup_{\xi\in Q_\ell^*}\fint_{Q_{\ell}^*}\frac{\|\gamma-g\|_{X\bigcap \mathcal{Y}(\xi')}}{|\xi-\xi'|^2}d\xi'.
\end{align*}

We now turn to the last estimate. Obviously,
\[
P^-_{\gamma}\gamma P^-_{\gamma}=P^-_{\gamma}(P^+_{g}-P^+_{\gamma})\gamma (P^+_g-P^+_{\gamma})P^-_\gamma.
\]
This implies
\begin{align*}
 \|P^-_{\gamma,\xi}\gamma_\xi P^-_{\gamma,\xi}\|_{X_c(\xi)}&\leq \frac{1+\kappa(\alpha,c)}{1-\kappa(\alpha,c)}\||\D_\xi|^{1/2}(P^+_{g,\xi}-P^+_{\gamma,\xi})\|_{\mathcal{B}(L^2_\xi)}^2\|\gamma_\xi\|_{\mathfrak{S}_1(\xi)}\\
    &\leq \frac{ C_{EE}^2}{2(1-\kappa)^{2}\lambda_0}\alpha_c^2\|g-\gamma\|_{X\bigcap \mathcal{Y}(\xi)}^2\|\gamma_\xi\|_{\mathfrak{S}_1(\xi)}.
\end{align*}
Thus,
\[
\|P^-_{\gamma}\gamma P^-_{\gamma}\|_{X_c}\leq \frac{ C_{EE}^2}{2(1-\kappa)^{2}\lambda_0}\alpha_c^2\fint_{Q_{\ell}^*}\|g-\gamma\|_{X\bigcap \mathcal{Y}(\xi)}^2 \|\gamma_\xi\|_{\mathfrak{S}_1(\xi)}d\xi.
\]
This ends the proof.
\end{proof}
Inserting Lemma \ref{lem:T-I} into \eqref{eq:6.3}, we get immediately
\begin{align*}
    |E(\gamma)-(\mathcal{E}(\gamma)-\epsilon_P\widetilde{\Tr}_{L^2}(\gamma))|\leq (2+\frac{\epsilon_P}{c^2})\alpha_c^2 \mathcal{N}_\gamma(\gamma-g),
\end{align*}
where $\mathcal{N}_\gamma(h)$ is given by \eqref{eq:error}, i.e.,
\begin{align*}
\mathcal{N}_\gamma(h)&=\frac{ C_{EE}^2}{2(1-\kappa)^{2}\lambda_0}\fint_{Q_{\ell}^*}\|h\|_{X\bigcap \mathcal{Y}(\xi)}^2 \|\gamma_\xi\|_{\mathfrak{S}_1(\xi)}d\xi\\
 &+(q\alpha_c^2+R\alpha_c^2+\alpha_c) \frac{10C_{EE}^4}{(1-\kappa)^4\lambda_0^{5/2}(1-L)^2} \\
  &\qquad\times\left(\frac{1}{c}\fint_{Q_{\ell}^*}\|h\|_{X\bigcap \mathcal{Y}(\xi)}\|\gamma_\xi|\D_\xi|^{1/2}\|_{\mathfrak{S}_1(\xi)}d\xi
+\frac{1}{c}\sup_{\xi\in Q_\ell^*}\fint_{Q_{\ell}^*}\frac{\|h\|_{X\bigcap \mathcal{Y}(\xi')}}{|\xi-\xi'|^2}d\xi'\right)^2.
\end{align*}

\subsection{Proof of Lemma \ref{lem:t-T}}\label{sec:6.3}
To complete the proof of Proposition \ref{main theorem}, it remains to prove Lemma \ref{lem:t-T}. Before going further, we need the following. 
\begin{proposition}
Let $\kappa<1$. For any $\gamma,\gamma'\in \Gamma_q$ and $h\in X$,
\begin{align}\label{eq:dP}
    P_{\gamma}^+ (dP^+_{\gamma}h) P^+_{\gamma}=0
\end{align}
where $dP^+_{\gamma}h$ is the Gateaux derivative of $P^+_{\gamma}$ at $\gamma\in\Gamma_q$ in the direction $h$. In addition, we have
\begin{align}\label{eq:5.11}
    \||\D|^{1/2}[P^+_{\gamma}-P^+_{\gamma'}-dP^+_{\gamma'}(\gamma-\gamma')]\|_Y\leq  \frac{C_{EE}^2}{2(1-\kappa)^{1/2}(1+C_{\mathcal{Y}})^2\lambda_0^{3/2}} \frac{\alpha_c^2}{c^3}\|\gamma-\gamma'\|_{X_c\bigcap \mathcal{Y}_c}^2.
\end{align}
\end{proposition}
\begin{proof}
The proof of this lemma is essentially the same as in \cite{meng}. As $P^+_{\gamma+th}=(P^+_{\gamma+th})^2$ for any $h\in X$, for any $h\in X$ we have
\[
dP^+_{\gamma}h=P^+_{\gamma}(dP^+_{\gamma}h)+(dP^+_{\gamma}h) P^+_{\gamma}.
\]
Thus,
\[
P^+_{\gamma}(dP^+_{\gamma}h)P^+_{\gamma}=2P^+_{\gamma}(dP^+_{\gamma}h)P^+_{\gamma},
\]
hence \eqref{eq:dP}. Recall that
\[
P^+_{\gamma}-P^+_{\gamma'}=\frac{\alpha}{2\pi}\int_{-\infty}^{+\infty}(\D_{\gamma}-iz)^{-1}V_{\gamma'-\gamma}(\D_{\gamma'}-iz)^{-1}dz,
\]
and
\[
dP^+_{\gamma}(\gamma-\gamma')=\frac{\alpha}{2\pi}\int_{-\infty}^{+\infty}(\D_{\gamma}-iz)^{-1}V_{\gamma'-\gamma}(\D_{\gamma}-iz)^{-1}dz.
\]
Thus,
\[
P^+_{\gamma}-P^+_{\gamma'}-dP^+_{\gamma}(\gamma-\gamma')=-\frac{\alpha^2}{2\pi}\int_{-\infty}^{+\infty}(\D_{\gamma}-iz)^{-1}V_{\gamma'-\gamma}(\D_{\gamma'}-iz)^{-1}V_{\gamma'-\gamma}(\D_{\gamma}-iz)^{-1}dz.
\]
Analogous to \eqref{eq:P-P}, using \eqref{3.eq:5.1'}, \eqref{eq:D-D}, \eqref{eq:arctan} and \eqref{eq:A.9} again, for any $\phi_\xi,\psi_\xi\in L^2_\xi$,
\begin{align*}
    \MoveEqLeft \left|\left(\phi_\xi,|\D_\xi|^{1/2}[P^+_{\gamma,\xi}-P^+_{\gamma',\xi}-(dP^+_{\gamma}(\gamma-\gamma'))_{\xi}]\psi_\xi\right)\right|\\
    &\leq \frac{\alpha^2}{2\pi}\|V_{\gamma'-\gamma}\|^2_Y\||\D_{\gamma'}|^{-1}\|_Y\\
    &\times\left(\int_{-\infty}^{+\infty}\|(\D_{\gamma,\xi}-iz)^{-1}|\D_{\xi}|^{1/2}\phi_\xi\|_{L^2_\xi}\right)^{1/2}\left(\int_{-\infty}^{+\infty}\|(\D_{\gamma,\xi}-iz)^{-1}\psi_\xi\|_{L^2_\xi}\right)^{1/2}\\
    &\leq \frac{C_{EE}^2}{2(1-\kappa)^{1/2}(1+C_{\mathcal{Y}})^2\lambda_0^{3/2}} \frac{\alpha_c^2}{c^3}\|\gamma-\gamma'\|_{X_c\bigcap \mathcal{Y}_c}^2\|\phi_\xi\|_{L^2_\xi}\|\psi_\xi\|_{L^2_\xi}.
\end{align*}
This gives \eqref{eq:5.11}. Hence the lemma.
\end{proof}

\begin{proof}[Proof of Lemma \ref{lem:t-T}]
We prove \eqref{theta-T+} first. Indeed, it suffices to prove
\begin{align*}
    \MoveEqLeft \|P^{+}_\gamma(T^n(\gamma)-T^{n-1}(\gamma))P^{+}_\gamma\|_{X_c\bigcap \mathcal{Y}_c}\\
&\leq C_{\kappa,L} (1-L) R \frac{\alpha_c^2}{c^2}\|T(\gamma)-\gamma\|_{X_c\bigcap \mathcal{Y}_c}\|T^{n-1}(\gamma)-T^{n-2}(\gamma)\|_{X_c\bigcap \mathcal{Y}_c}.
\end{align*}
Then by Lemma \ref{lem:retra}, 
\begin{align*}
 \MoveEqLeft \|P^{+}_\gamma(\theta(\gamma)-T(\gamma))P^{+}_\gamma\|_{X_c\bigcap \mathcal{Y}_c}\leq \sum_{n=2}^{+\infty}\|P^{+}_\gamma(T^n(\gamma)-T^{n-1}(\gamma))P^{+}_\gamma\|_{X_c\bigcap \mathcal{Y}_c}\\
  &\leq C_{\kappa,L}(1-L) R\frac{\alpha_c^2}{c^2}\|T(\gamma)-\gamma\|_{X_c\bigcap \mathcal{Y}_c}\sum_{n=2}^{+\infty}\|T^{n-1}(\gamma)-T^{n-2}(\gamma)\|_{X_c\bigcap \mathcal{Y}_c}\\
  &\leq C_{\kappa,L} R\frac{\alpha_c^2}{c^2}\|T(\gamma)-\gamma\|_{X_c\bigcap \mathcal{Y}_c}^2.
\end{align*}
Let $\gamma_n=T^n(\gamma)$ and $\gamma_0=\gamma$. Then for $n\geq 2$, $\gamma_n=P^+_{\gamma_{n-1}}\gamma_{n-1}P^+_{\gamma_{n-1}}$ and $\gamma_{n-1}=P^+_{\gamma_{n-2}}\gamma_{n-1}P^+_{\gamma_{n-2}}$. Hence, for $n\geq 2$
\begin{align}\label{eq:decom'}
 P^+_\gamma(\gamma_n-\gamma_{n-1})P^+_\gamma&=P^+_\gamma(P^+_{\gamma_{n-1}}-P^+_{\gamma_{n-2}})P^+_{\gamma_{n-2}}\gamma_{n-1}P^+_{\gamma_{n-1}}P^+_\gamma\notag\\
 &\qquad+P^+_\gamma\gamma_{n-1}P^+_{\gamma_{n-2}}(P^+_{\gamma_{n-1}}-P^+_{\gamma_{n-2}})P^+_\gamma.
\end{align}
We only need to consider the first term on the right-hand side, the second term can be treated in the same manner. According to \eqref{eq:dP}, we have
    \begin{align}\label{eq:decom}
\MoveEqLeft P^+_\gamma(P^+_{\gamma_{n-1}}-P^+_{\gamma_{n-2}})P^+_{\gamma_{n-2}}\gamma_{n-1}P^+_{\gamma_{n-1}}P^+_\gamma\notag\\
    &=P^+_{\gamma_{n-2}}(P^+_{\gamma_{n-1}}-P^+_{\gamma_{n-2}}-dP^+_{\gamma_{n-2}}(\gamma_{n-1}-\gamma_{n-2}))P^+_{\gamma_{n-2}}\gamma_{n-1}P^+_{\gamma_{n-1}}P^+_\gamma\notag\\
    &\quad+(P^+_\gamma-P^+_{\gamma_{n-2}})(P^+_{\gamma_{n-1}}-P^+_{\gamma_{n-2}})P^+_{\gamma_{n-2}}\gamma_{n-1}P^+_{\gamma_{n-1}}P^+_\gamma.
\end{align}
Thus, according to \eqref{eq:5.11} and \eqref{eq:DPD}, for the first term on the right-hand side
\begin{align*}
    \MoveEqLeft \|P^+_{\gamma_{n-2}}(P^+_{\gamma_{n-1}}-P^+_{\gamma_{n-2}}-dP^+_{\gamma_{n-2}}(\gamma_{n-1}-\gamma_{n-2}))P^+_{\gamma_{n-2}}\gamma_{n-1}P^+_{\gamma_{n-1}}P^+_\gamma\|_{X_c\bigcap \mathcal{Y}_c}\\
    &\leq \frac{(1+\kappa)^{1/2}}{(1-\kappa)^{1/2}}\||\D|^{1/2}(P^+_{\gamma_{n-1}}-P^+_{\gamma_{n-2}}-dP^+_{\gamma_{n-2}}(\gamma_{n-1}-\gamma_{n-2}))\|_Y\\
    &\quad \times\max\{\|\gamma_{n-1}P^+_{\gamma_{n-1}}P^+_{\gamma}|\D|^{1/2}\|_{\mathfrak{S}_1}+\|\gamma_{n-1}\|_{\mathcal{Y}}\}\\
    &\leq \frac{C_{EE}^2(1+\kappa)^{3/2} }{2(1-\kappa)^2(1+C_{\mathcal{Y}})^2\lambda_0^{3/2}} \frac{\alpha_c^2}{c^3}\|\gamma_{n-1}-\gamma_{n-2}\|_{X_c\bigcap \mathcal{Y}_c}^2\max\{\|\gamma_{n-1}|\D|^{1/2}\|_{\mathfrak{S}_{1,1}}+\|\gamma_{n-1}\|_{\mathcal{Y}}\}.
\end{align*}
As $\gamma\in \mathcal{U}_R$ and according to Lemma \ref{lem:retra} for any $n\geq 2$,
\[
\|\gamma_{n-1}-\gamma_{n-2}\|_{X_c\bigcap \mathcal{Y}_c}\leq \|\gamma-T(\gamma)\|_{X_c\bigcap \mathcal{Y}_c},\quad \frac{1}{c}\max\{\|\gamma_{n-1}|\D|^{1/2}\|_{\mathfrak{S}_{1,1}},\|\gamma_{n-1}\|_{\mathcal{Y}}\}\leq R.
\]
Then as $\kappa\leq 1$ and $C_{\mathcal{Y}}\geq 0$, we have
\begin{align*}
    \MoveEqLeft \|P^+_{\gamma_{n-2}}(P^+_{\gamma_{n-1}}-P^+_{\gamma_{n-2}}-dP^+_{\gamma_{n-2}}(\gamma_{n-1}-\gamma_{n-2}))P^+_{\gamma_{n-2}}\gamma_{n-1}P^+_{\gamma_{n-1}}P^+_\gamma\|_{X_c\bigcap \mathcal{Y}_c}\\
    &\qquad\qquad\qquad\qquad\leq  C_{\kappa,L}' R\frac{\alpha_c^2}{c^2}\|T(\gamma)-\gamma\|_{X_c\bigcap \mathcal{Y}_c}\|\gamma_{n-1}-\gamma_{n-2}\|_{X_c\bigcap \mathcal{Y}_c},
\end{align*}
with $C_{\kappa,L}':=\frac{2C_{EE}^2}{(1-\kappa)^2\lambda_0^{3/2}}$. 

Now we consider the second term on the right-hand side of \eqref{eq:decom}. By \eqref{eq:3.3} we have $\|\gamma_{n}-\gamma\|_{X_c\bigcap \mathcal{Y}_c}\leq \frac{1}{1-L}\|T(\gamma)-\gamma\|_{X_c\bigcap \mathcal{Y}_c}$. Thus, by \eqref{eq:P-P},
\begin{align*}
\MoveEqLeft    \|(P^+_\gamma-P^+_{\gamma_{n-2}})(P^+_{\gamma_{n-1}}-P^+_{\gamma_{n-2}})P^+_{\gamma_{n-2}}\gamma_{n-1}P^+_{\gamma_{n-1}}P^+_{\gamma}\|_{X_c\bigcap \mathcal{Y}_c}\\
&\leq\frac{(1+\kappa)}{(1-\kappa)}\||\D|^{1/2}(P^+_\gamma-P^+_{\gamma_{n-2}})\|_Y\||\D|^{-1/2}\|_Y\||\D|^{1/2}(P^+_{\gamma_{n-1}}-P^+_{\gamma_{n-2}})\|_Y\\
&\quad\times\max\{\|\gamma_{n-1}|\D|^{1/2}\|_{\mathfrak{S}_1},\|\gamma_{n-1}\|_{\mathcal{Y}}\}\\
&\leq \frac{C_{EE}^2(1+\kappa)}{4(1-\kappa)^2(1+C_{\mathcal{Y}})^2\lambda_0^{3/2}} R \frac{\alpha_c^2}{c^2}\|\gamma_{n-2}-\gamma\|_{X_c\bigcap \mathcal{Y}_c}\|\gamma_{n-1}-\gamma_{n-2}\|_{X_c\bigcap \mathcal{Y}_c}\\
&\leq C_{\kappa,L}'' R\frac{\alpha_c^2}{c^2}\|T(\gamma)-\gamma\|_{X_c\bigcap \mathcal{Y}_c}\|\gamma_{n-1}-\gamma_{n-2}\|_{X_c\bigcap \mathcal{Y}_c},
\end{align*}
with $C_{\kappa,L}'':=\frac{C_{EE}^2}{2(1-\kappa)^2\lambda_0^{3/2}(1-L)}$. Thus by \eqref{eq:decom},
\begin{align*}
   \MoveEqLeft \|P^+_\gamma(P^+_{\gamma_{n-1}}-P^+_{\gamma_{n-2}})P^+_{\gamma_{n-2}}\gamma_{n-1}P^+_{\gamma_{n-1}}P^+_\gamma\|_{X_c\bigcap\mathcal{Y}_c}\\
    &\leq (C_{\kappa,L}'+C_{\kappa,L}'')R\frac{\alpha_c^2}{c^2}\|T(\gamma)-\gamma\|_{X_c\bigcap \mathcal{Y}_c}\|T^{n-1}(\gamma)-T^{n-2}(\gamma)\|_{X_c\bigcap \mathcal{Y}_c}.
\end{align*}
The second terms on the right-hand side of \eqref{eq:decom'} can be treated analogously, thus
\[
\|P^+_\gamma(\gamma_n-\gamma_{n-1})P^+_\gamma\|_{X_c\bigcap \mathcal{Y}_c}\leq C_{\kappa,L}(1-L)R\frac{\alpha_c^2}{c^2}\|T(\gamma)-\gamma\|_{X_c\bigcap \mathcal{Y}_c}\|T^{n-1}(\gamma)-T^{n-2}(\gamma)\|_{X_c\bigcap \mathcal{Y}_c}.
\]
where $C_{\kappa,L}:= \frac{5C_{EE}^2}{(1-\kappa)^2\lambda_0^{3/2}(1-L)^2}\geq 2 (1-L)^{-1}(C_{\kappa,L}'+C_{\kappa,L}'')$. Hence \eqref{theta-T+}.

\medskip

Finally, we consider the term $P^-_\gamma\theta(\gamma)P^-_\gamma$. As $\theta(\gamma)=P^+_{\theta(\gamma)}\theta(\gamma)P^+_{\theta(\gamma)}$, we have
\begin{align*}
    P^-_\gamma\theta(\gamma)P^-_\gamma= P^-_{\gamma}(P^+_{\theta(\gamma)}-P^+_{\gamma})\theta(\gamma)(P^+_{\theta(\gamma)}-P^+_{\gamma})P^-_{\gamma},
\end{align*}
from which we deduce
\begin{align*}
\|P^-_\gamma\theta(\gamma)P^-_\gamma\|_{X_c\bigcap \mathcal{Y}_c}&\leq \frac{1+\kappa}{1-\kappa}\||\D|^{1/2}(P^+_{\theta(\gamma)}-P^+_{\gamma})\|_Y^2\max\{\|\theta(\gamma)\|_{\mathfrak{S}_{1,1}},\|\theta(\gamma)\|_{\mathcal{Y}}\}\\
    &\leq \frac{C_{EE}^2(1+\kappa)}{4(1-\kappa)^2(1+C_\mathcal{Y})^2\lambda_0 (1-L)^2}(q+C_{\mathcal{Y}})\frac{\alpha_c^2}{c^2}\|T(\gamma)-\gamma\|_{X_c\bigcap \mathcal{Y}_c}^2\\
    &\leq  C_{\kappa,L}q\frac{\alpha_c^2}{c^2}\|T(\gamma)-\gamma\|_{X_c\bigcap \mathcal{Y}_c}^2.
\end{align*}
In the last inequality, we use $(1+C_\mathcal{Y})^{-2}(q+C_{\mathcal{Y}})\leq q(1+C_{\mathcal{Y}})^{-1}\leq 1$. This ends the proof.
\end{proof}

\appendix

\section{Some useful inequalities}
In this appendix, we recall some useful results proved in \cite{crystals}, and we adapt these results for the new norm $\mathcal{Y}$. 
\begin{lemma}[Hardy-type inequalities]\label{lem:hardy}
There exist positive constants $C_G:=C_G(\ell)\geq 1$, $C_W:=C_W(\ell)\geq 1$, $C_{EE}:=C_{EE}(\ell)\geq 1$ and $C_{EE}':=C_{EE}'(\ell)$ that only depend on $\ell$ and such that for any $\xi\in Q_\ell^*$,
\begin{align}\label{eq:W-Y}
    \|W_{\gamma,\xi}\|_{\mathcal{B}(L^2_\xi)}\leq C_W\|\gamma\|_{X\bigcap \mathcal{Y}(\xi)}\leq \frac{C_W}{c}\|\gamma\|_{X_c\bigcap \mathcal{Y}_c(\xi)},
\end{align}
\begin{equation}\label{3.eq:4.3}
\begin{aligned}
\MoveEqLeft \|G_{\ell}\psi_\xi\|_{L^2_\xi}\leq  C_G\,\|(1-\Delta_\xi)^{1/2}\psi_\xi\|_{L^2_\xi}\leq \frac{C_G}{c}\||\D_\xi|\psi_\xi\|_{L^2_\xi},
\end{aligned}
\end{equation}
\begin{align}\label{3.eq:5.1''}
    \|V_{\gamma}\|_{Y} \leq C_{EE}\,\|\gamma\|_{X\bigcap Y} \leq \frac{C_{EE}}{c}\,\|\gamma\|_{X_c\bigcap Y_c},
\end{align}
\begin{align}\label{3.eq:5.1'}
\|V_{\gamma}\|_{Y} \leq \frac{C_{EE}}{(1+C_\mathcal{Y})}\,\|\gamma\|_{X\bigcap \mathcal{Y}}\leq \frac{C_{EE}}{(1+C_\mathcal{Y})\,c}\,\|\gamma\|_{X_c\bigcap \mathcal{Y}_c},
\end{align}
\begin{align}\label{3.eq:5.1}
\|V_{\gamma}\,(1-\Delta)^{-1/2}\|_{Y} \leq C_{EE}\,\|\gamma\|_{\mathfrak{S}_{1,1}\bigcap Y},
\end{align}
and
\begin{align}\label{3.eq:5.1'''}
     -C_{EE}'\,\|\gamma\|_{\mathfrak{S}_{1,1}\bigcap Y}\|\psi_\xi\|_{L^2_\xi}^2\leq \left(\psi_\xi,V_{\gamma,\xi}\psi_\xi\right)_{L^2_\xi}.
\end{align}
\end{lemma}
\begin{proof}
These estimates can be found in \cite[Lemma 4.5 and Lemma 4.7]{crystals}. The estimate \eqref{eq:W-Y} is slightly different from the one in \cite{crystals} because of the change of the functional space ($\mathcal{Y}$ instead of $Y$). According to \cite[Eq. (B.13)-(B.14)]{crystals}, we know that, for any $\xi\in Q_\ell^*$,
\begin{align*}
\left\|\fint_{Q_\ell^*}d\xi'\int_{Q_\ell}W_{<2,\ell}^\infty(\xi'-\xi,x-y)\gamma_{\xi'}(x,y)\psi_\xi(y)\,dy\right\|_{L^2_\xi}&\leq \frac{\|\psi_\xi\|_{L^2_\xi}}{2\pi^2}\,\int_{3Q_\ell^*}\frac{\|\gamma_{\xi'}\|_{\mathcal{B}(L^2_{\xi'})}}{|\xi'-\xi|^2}\, d\xi'\\
   &\leq \frac{27}{2\pi^2}\,\|\gamma\|_{\mathcal{Y}(\xi)}\|\,\psi_\xi\|_{L^2_\xi},
\end{align*}
since $\|\gamma_{\xi'}\|_{\mathcal{B}(L^2_{\xi'})}=\|\gamma_{\xi'+p}\|_{\mathcal{B}(L^2_{\xi'+p})}$ with $p\in \frac{2\pi}{\ell}\mathbb{Z}^3$. Repeating the proof of \cite[Eq. (4.8)]{crystals}, the estimate \eqref{eq:W-Y} follows. Then, \eqref{3.eq:5.1'} is a modification of \eqref{3.eq:5.1''} where we replace the estimates on the exchange term $W_{\gamma}$ by \eqref{eq:W-Y}.
\end{proof}
\begin{lemma}\cite[Lemma 4.9,  Eq.(5.13) and Lemma 4.10]{crystals}\label{lem:ope}
Let $\gamma\in \Gamma_{\leq q}$.
\begin{enumerate}
    \item  If $\kappa(\alpha,c)<1$, then
    \begin{align*}
        (1-\kappa(\alpha,c))^2|\D|^2\leq |\D_{\gamma}|^2\leq  (1+\kappa(\alpha,c))^{2}\vert
        \D\vert^2,
    \end{align*}
in the sense of operators. Consequently,
\begin{align}\label{eq:D-D}
        (1-\kappa(\alpha,c))|\D|\leq |\D_{\gamma}|\leq  (1+\kappa(\alpha,c))|\D|.
    \end{align}
\item Let $\gamma\in \Gamma_q$ and $\kappa(\alpha,c)<1$. Then
    \begin{align}\label{eq:DPD}
        \||\D|^{1/2} P^\pm_{\gamma}|\D|^{-1/2}\|_{Y}\leq
\frac{(1+\kappa(\alpha,c))^{1/2}}{(1-\kappa(\alpha,c))^{1/2}}.
\end{align}
\item Let $\gamma\in \Gamma_q$ and $\kappa(\alpha,c)<1$, then
\begin{align}\label{eq:A.9}
    \inf|\sigma(\D_{\gamma})|\geq c^2 \lambda_0(\alpha,c)\geq c^2(1-\kappa(\alpha,c)),
\end{align}
with $\lambda_0(\alpha,c):=1- \max\big\{C_H z_c+C_{EE}' \,\alpha_c q, \frac{C_0}{\ell}z_c+C_{EE}\, \alpha_c q\big\}$.
\end{enumerate}
\end{lemma}

\section{Boundedness of eigenfunctions of \texorpdfstring{$\D_{\gamma,\xi}$}{}}
We start with the following. 
\begin{lemma}[Properties of  positive eigenvalues of $\D_{\gamma,\xi}$]\label{lem:prop-spec}
Assume that $c>C_Gz$ and $\kappa(\alpha,c)<1$. Let  $\gamma\in \Gamma_{\leq q}$. For $k\geq 1$, we denote by $\lambda_k(\xi)$ the  $k$-th positive eigenvalue (counted with multiplicity)  of the mean-field operator $\D_{\gamma,\xi}$. Then, there exist positive  constants $c^*(k)$ and $c_*(k)$ independent of $\xi$ and $\Sigma(k)$ independent of $\alpha$, $c$ and $\xi$, such that 
\begin{align}\label{eq:ck-ck}
    0\leq c_*(k)\leq \lambda_k(\xi)\leq c^*(k)\leq c^2+\Sigma(k),
\end{align}
for every $\xi \in Q_{\ell}^*$, with $c_*(k)\to+\infty$ when $k\to+\infty$. The interval $[c_*(k),c^*(k)]$ is  independent of $\gamma$ in $\Gamma_{\leq q}$. In particular, there exists an integer $M>0$ independent of $\xi$ such that there are at most $q+M$ eigenvalues in $[0,c^*(q+1)]$. 
\end{lemma}
\begin{proof}[Sketch of proof]
The proof is in the very same spirit as for \cite[Lemma 4.11]{crystals}. Additionally, we show here that there exists $\Sigma(k)$ independent of $\alpha$, $c$ and $\xi$ such that \eqref{eq:ck-ck} holds. (Note that compared with \cite{crystals}, we use different definitions for $c_*$ and $c^*$.) We modify the proof of Lemma 4.8 in \cite{crystals} by using the Hardy-type inequalities associated with the operator $(1-\Delta_\xi)^{1/2}$ instead of $|\D_\xi|$ in Lemma \ref{lem:hardy}.  From \eqref{3.eq:4.3} and \eqref{3.eq:5.1}, we get
\begin{align*}
    \MoveEqLeft\lambda_k(\xi)\geq \sigma_k(|\D_\xi|-(C_Gz+\alpha C_{EE}q)(1-\Delta_\xi)^{1/2}),
\end{align*}
where $\sigma_k(A)$ is the $k$-th positive eigenvalue (counted with multiplicity) of the operator $A$. Hence,
\[
c_*(k):=\inf_{\xi}\sigma_k((\D_\xi-(C_Gz+\alpha C_{EE}q)(1-\Delta_\xi)^{1/2})).
\]
On the other hand,
\begin{align*}
    \lambda_k(\xi)&\leq \sigma_k((\D_\xi+(C_Gz+\alpha C_{EE}q)(1-\Delta_\xi)^{1/2}))\\
    &\leq d^+_k(\xi)+(C_Gz+\alpha C_{EE}q)\sigma_k((1-\Delta_\xi)^{1/2}).
\end{align*}
where the second inequality holds by using  the min–max Courant--Fisher formula. We may choose 
\[
c^*(k):=\sup_{\xi}\left(d^+_k(\xi)+(C_Gz+\alpha C_{EE}q)\sigma_k((1-\Delta_\xi)^{1/2})\right).
\]
Notice that there is $m(k)\in\mathbb{N}^3$ such that
\[
d^+_k(\xi)=c^2\sqrt{1+\frac{1}{c^2}\left|\xi+\frac{2\pi}{\ell}m(k)\right|^2}\leq c^2+\frac{1}{2}\left|\xi+\frac{2\pi}{\ell}m(k)\right|^2
\]
and
\begin{align*}
    \sigma_k((1-\Delta_\xi)^{1/2})\leq \left|\xi+\frac{2\pi}{\ell}m(k)\right|.
\end{align*}
Therefore, since $\alpha \leq 1$, \eqref{eq:ck-ck} holds with 
\[\Sigma(k)=\frac{2\pi^2(1+|m(k)|)^2 }{\ell^2}+\left(C_Gz+C_{EE}q\right)\frac{2\pi(1+|m(k)|) }{\ell}.
\]
This proves the property of $\Sigma(k)$.
\end{proof}
Then we have the following a priori estimate on $H^1$ norms of eigenfunctions. 
\begin{lemma}\label{lem:unibound}
Let $\gamma\in \Gamma_{\leq q}$ and let $\alpha$ and $c$ be chosen such that $\kappa(\alpha,c)\leq 1-C_1$ for some $0<C_1<1$ independent of $\alpha$ and $c$.  Let also $\psi_\xi$ be a normalized eigenfunction of the operator $\D_{\gamma,\xi}$ with eigenvalue  $\lambda(\xi) \in (0,c^2+\Sigma(q+1)]$  where $\Sigma(q+1)$ is given in Lemma \ref{lem:prop-spec}. Then,  there exists a positive constant $K$ that is independent of $\xi$, $\alpha$ and $c$ such that
\[
\|\psi_\xi\|_{L^\infty(Q_{\ell}^*;H^1_\xi(Q_\ell))}\leq K.
\]
Furthermore, if $\gamma'\in \Gamma_q$ satisfies  $0\leq \gamma' \leq \mathbbm{1}_{(0,c^2+\Sigma(q+1))}(\D_{\gamma})$, then
\[
\|\gamma'\|_X\leq K^2q.
\]
\end{lemma}
\begin{proof}
As $\D_{\gamma,\xi}\psi_{\xi}=\lambda(\xi) \psi_{\xi}$, we have
\[
\|\D_\xi \psi_{\xi}\|_{L^2_\xi}=\|(\lambda(\xi)+zG_\ell-\alpha W_{\gamma,\xi})\psi_\xi\|_{L^2_\xi}.
\]
According to the Hardy inequality and since $|\lambda(\xi)|\leq c^2+\Sigma(q+1)$,
\begin{align*}
   c^4\|\psi_\xi\|^2_{L^2_\xi}+c^2\|\nabla_\xi \psi_\xi\|^2_{L^2_\xi}&\leq (c^4+2c^2\Sigma(q+1)+\Sigma(q+1)^2)\|\psi_\xi\|^2_{L^2_\xi}\\
    &\quad +2( C_G z+\alpha C_{EE}q)(c^2+\Sigma(q+1))\|(1-\Delta_\xi)^{1/2} \psi_\xi\|_{L^2_\xi}\|\psi_\xi\|_{L^2_\xi} \\
    &\quad+c^2\kappa(\alpha,c)^2\|(1-\Delta_\xi)^{1/2}\psi_\xi\|^2_{L^2_\xi}.
\end{align*}
 As $\|\psi\|_{L^2_\xi}=1$, we obtain 
 \begin{align*}
\|\nabla_\xi \psi_\xi\|^2_{L^2_\xi}&\leq  \big(2\Sigma(q+1)+\frac{1}{c^2}\Sigma(q+1)^2\big)+ \kappa(\alpha,c)^2\|(1-\Delta_\xi)^{1/2}\psi_\xi\|^2_{L^2_\xi}\\
    &\quad+2( C_G z+\alpha C_{EE}q)\Big(1+\frac{1}{c^2}\Sigma(q+1)\Big)\|(1-\Delta_\xi)^{1/2} \psi_\xi\|_{L^2_\xi}.
 \end{align*}
According to our assumption, $\kappa(\alpha,c)\leq 1-C_1<1$. Then using the Young inequality for the last term on the right-hand side, we infer that for some $K$ large enough, independent of $c$ and $\xi$,
 \begin{align*} \|(1-\Delta_\xi)^{1/2} \psi_\xi\|^2_{L^2_\xi} \leq K^2.
 \end{align*}
For the second estimate, we express $\gamma'$ as 
\[
\gamma'_\xi=\sum_{i=1}^{+\infty} \mu_i(\xi)\left|\psi_i(\xi)\right>\left<\psi_i(\xi)\right|
\]
where $0\leq \mu_i(\xi)\leq 1$, $\sum_{i=1}^{+\infty}\fint_{Q_\ell^*}\mu_i(\xi)\,d\xi=q$ and $\psi_i(\xi)$ is the normalized eigenfunction of $\D_{\gamma,\xi}$ with eigenvalue $\lambda_i(\xi) \in (0; c^2+\Sigma(q+1))$. Thus,
\[
\|\gamma'\|_{X}\leq \sum_{i=1}^{+\infty}\fint_{Q_\ell^*}\mu_i(\xi)\|\psi_i(\xi)\|_{H^1_\xi}^2d\xi\leq K^2q.
\]
This concludes the proof.
\end{proof}


\noindent\textbf{Acknowledgments.} L. M. acknowledges support from the European Research Council (ERC) under the European Union's Horizon 2020 research and innovation program (grant agreement No. 810367).


\begin{refcontext}[sorting=nyt]
\printbibliography[heading=bibintoc, title={Bibliography}]
\end{refcontext}

\end{document}